\documentclass[10pt]{article}
\usepackage{amssymb}
\usepackage{amsthm}
\usepackage{epsf,epsfig,amsfonts,graphicx,color}
\usepackage{a4wide}
\usepackage{epsf,epsfig,amsfonts}

\usepackage{amsfonts}
\usepackage{amsmath,amssymb,amsthm}
\usepackage{url}\usepackage{amsmath,amssymb,amsthm}
\usepackage{url}

\newtheorem{Th}{Theorem}
\newtheorem{Prop}{Proposition}
\newtheorem{Lemma}{Lemma}
\newtheorem{Cor}{Corollary}

\theoremstyle{remark}
\newtheorem{Ex}{Example}
\newtheorem{Rem}{Remark}

\newcommand{\const}{\mbox{\rm const}}

\newcommand{\trace}{\mbox{\rm trace}}

\newcommand{\eusl}{\mathfrak{sl}}

\newcommand{\bbR}{\mathbb{R}}

\newcommand{\weg}[1]{}
\newcommand{\sign}{\mbox{\rm sign}}\newcommand{\Proj}{\mbox{\rm Proj}}\newcommand{\Iso}{\mbox{\rm Iso}}
\def\<{\langle}
\def\>{\rangle}

\begin{document}
\title%[Pseudo-Riemannian metric admitting quadratic integrals]
{Pseudo-Riemannian metrics on closed surfaces whose geodesic flows admit nontrivial  integrals  quadratic in momenta, and proof of the  projective Obata conjecture for two-dimensional pseudo-Riemannian metrics}
\author{Vladimir S. Matveev\thanks{Institute of Mathematics, FSU Jena, 07737 Jena Germany,  vladimir.matveev@uni-jena.de}\ \thanks{Partially supported by DFG (SPP 1154 and GK 1523)  }}
\date{}
\maketitle

\begin{abstract} 
We describe all pseudo-Riemannian metrics on  closed surfaces whose geodesic flows  admit nontrivial  integrals quadratic in momenta. As an application, we  solve  the Beltrami problem on closed surfaces,  prove  the nonexistence of quadratically-superintegrable metrics of nonconstant curvature 
 on closed surfaces, and prove the two-dimensional pseudo-Riemannian 
 version of the projective Obata conjecture.  
\end{abstract}

\section{Introduction} 
\subsection{Definitions and  the  statement of the problem} 
Consider a  pseudo-Riemmanian metric $g=(g_{ij})$ on  a surface $M^2$.  A function $F:{T^*}M\to \mathbb{R}$ is called {\it an integral} of the geodesic flow of $g$, if $\{ H, F\}=0$, where $H:= \tfrac{1}{2} \sum_{i,j} g^{ij} p_ip_j:T^*M\to \mathbb{R}$ is the kinetic energy corresponding to the metric. Geometrically, the condition  $\{ H, F\}=0$  means 
that the function $F$ is constant on the trajectories  of the Hamiltonian system with the Hamiltonian $H$. We say that  the  integral $F$ is {\it quadratic in  momenta,} if in every local coordinate system $(x,y)$ on $M^2$  it has the form 
\begin{equation} \label{integral} 
a(x,y)p_x^2+ b(x,y)p_xp_y+ c(x,y)p_y^2  
\end{equation} in the canonical coordinates $(x,y,p_x, p_y)$ on $T^*M^2$. 
Geometrically, the  formula \eqref{integral} 
means that the restriction of the integral to every cotangent space $T^*_{(x,y)}M^2\equiv \mathbb{R}^2$  is a homogeneous quadratic function.  As  {\it trivial}  examples  of quadratic in momenta  integrals 
we consider those proportional to the Hamiltonian  $H$. 

Similarly, we say that the integral is {\it linear} in momenta, if  for  every local coordinate system $(x,y)$ on $M^2$  it has the form $\alpha(x,y) p_x + \beta(x,y) p_y$ in the canonical coordinates $(x,y,p_x, p_y)$ on $T^*M^2$; an integral linear in momenta is   {\it trivial}, if it is identically zero.

The importance of integrals quadratic in momenta  other than the Hamiltonian 
for studying the metric  was recognized long ago. Indeed,
it was Jacobi's realization that the geodesic flow of the ellipsoid
admitted such an `extra' quadratic integral that allowed him to
integrate the geodesics on the ellipsoid.

In the present paper we  solve (see Model Examples 1, 2, 3  and Theorems \ref{main1}, \ref{main3}, \ref{main4} below) the following problem: 

{\bf Problem.} {\it Find all metrics of signature $(+,-)$ on closed 2-dimensional manifolds  whose    geodesic flows  admit  nontrivial
integrals  quadratic in momenta.  } 

 Riemannian metrics whose geodesic flows  admit integrals quadratic in momenta are quite good studied. Indeed, local description of such  a  metric in a neighborhood of almost every point is known since Liouville. 
Moreover,  the Riemannian version (and, therefore, if the signature of $g$ is  (--,--)) of  the   problem above  was  solved. There exist two different approaches that lead to a solution:  one, which is based on  the  ideas  of  Kolokoltsov \cite{Kol},  was realized   in \cite{Kol,1,klein}, see  also \cite{BMF,BF}. 
Alternative  approach to the description of metrics whose geodesic flows admit  nontrivial integrals  quadratic   in momenta 
is due to Kiyohara \cite{Kio}, see also \cite{Igarashi,memo}.   Our solution uses main  ideas from both approaches.  

Metrics whose geodesic flows  admit integrals quadratic in momenta were studied in the framework of differential geometry (at least since Darboux \cite{Darboux}) and mathematical physics (at least since Birkhoff  \cite{3} and Whittaker \cite{Whi}). We give two 
 applications of our results in differential geometry and one application  in mathematical physics. In differential geometry, we use the connection between   integrals quadratic in momenta and geodesically equivalent metrics (we give the necessary definition in \S \ref{beltrami}) to solve the natural generalization of the Beltrami problem for closed manifolds, and to prove the two-dimensional pseudo-Riemannian 
 version of the projective Obata conjecture.  In mathematical physics, we prove that  all quadratically-superintegrable metrics on closed surfaces  (the necessary definition is in \S \ref{superi}) have constant curvature.  This generalizes  the result of \cite{Kio,klein} to the pseudo-Riemannian metrics.

\subsection{Metrics on the torus  whose geodesic flows admit nontrivial  integrals  quadratic in momenta}
Locally, pseudo-Riemannian metrics admitting integrals quadratic in momenta were described\footnote{As it mentioned in \cite{pucacco,BMP}, the essential part of the result appeared already in  Darboux  \cite[\S\S592--594,600--608]{Darboux}} in \cite[Theorem 1]{pucacco} and \cite[Theorem 1]{BMP}:

\begin{Th}[\cite{pucacco,BMP}] \label{main}   Suppose a Riemannian  or pseudo-Riemannian 
metric $g$  on a connected surface  ${M}^2$ admits an integral $F$ quadratic in momenta such that $F\ne \const \cdot H$ for all $\const \in \mathbb{R}$. Then, in a neighbourhood of almost every point there exist coordinates  $x,y$ such that  the metric and the integral are as in the following table:
\begin{center}\begin{tabular}{|c||c|c|c|}\hline & \textrm{Liouville case} & \textrm{Complex-Liouville case} & \textrm{Jordan-block case}\\ \hline \hline
$g$ & $(X(x)-Y(y))(dx^2 + \varepsilon dy^2)$ &  $\Im(h)dxdy$ & $\left(\widehat Y(y)+\tfrac{x}{2} Y'(y)\right)dxdy $
\\ \hline  $F$&$\tfrac{X(x)p_y^2 +\varepsilon Y(y)p_{x}^2}{X(x)-Y(y)} $&
 $p_x^2 - p_y^2   + 2\tfrac{\Re(h)}{\Im(h)}p_xp_y $ & $\varepsilon \left(p_x^2 -  \frac{Y(y)}{\widehat Y(y)+ \tfrac{x}{2} Y'(y) }p_xp_y\right)$  \\ \hline
\end{tabular}
\end{center}
where $\varepsilon = \pm1$, and  $\Re (h)$ and $\Im (h)$ are the real and imaginary parts of a
holomorphic function $h$ of the variable $z:= x+ i\cdot y$.
\end{Th}

\begin{Rem}  Within our paper, we understand ``{\it almost every}'' in the topological sense:
a condition is fulfilled at almost every point, if the set of the points where it is fulfilled is everywhere dense.  
\end{Rem} 

We see that the metric $g$ in the Jordan-block and Complex-Liouville cases always  has indefinite signature 
(+,--),   and the metric $g$ in the Liouville case has signature (+,--) if and only if $\varepsilon =-1$. 
The Liouville case with $\varepsilon =1$  was known to classics.

In Section \ref{puc}, we repeat  the proof of   Theorem \ref{main}, because we will need most techical details from   it in the proof of our main result, which is Theorem \ref{main1} below.   

Let us now discuss the case when $M^2$ is closed. 
First of all, because of Euler characteristic, a closed surface admitting a pseudo-Riemannian metric of indefinite signature   is homeomorphic to the torus or to the Klein bottle. Since a double cover of  the  Klein  bottle   is  the torus, and the geodesic flow of the  lift of a  metric  whose geodesic flow admits  an  integral quadratic in momenta  also admits an integral quadratic in momenta, the most  important  case is when   the  surface is the torus. In Model Example 1  below we describe a class of pseudo-Riemannian metrics on the torus such that their geodesic flows admit  nontrivial integrals  quadratic in momenta. Theorem \ref{main1} claims that every metric such that its  geodesic flow admits  a nontrivial integral quadratic in momenta  is isometric to one  from Model Example 1.

{\bf Model Example 1.} We consider $\mathbb{R}^2 $ with the standard coordinates $(x,y)$,  two linearly  independent vectors $\xi=(\xi_1,\xi_2)$ and $\nu=(\nu_1, \nu_2)$, and two nonconstant functions $X$ and $Y$ of one variable (it is convenient to think that the variable of $X$ is $x$ and the variable of $Y$ is $y$)
 such that  \begin{itemize} 
 \item[(a)] for all  $(x,y)\in \mathbb{R}^2$ we have $X(x)\ne Y(y)$, and 
   \item[(b)] for every $(x,y)\in \mathbb{R}^2$, $ X(x+\xi_1)= X(x+ \nu_1)= X(x)$ and $ Y(y+\xi_2)= Y(y+ \nu_2)= Y(y)$. \end{itemize}
 Next, consider the metrics $(X(x)- Y(y))(dx^2 +\varepsilon  dy^2)$  on  $\mathbb{R}^2$, where $\varepsilon =\pm1$,  and 
 the action of the lattice  $G:=\{k \cdot \xi  + m \cdot \nu \mid k, m\in \mathbb{Z}\}$ on $\mathbb{R}^2$. The action is free, discrete and 
 preserves the metric and the quadratic integral  $ \tfrac{X(x)p_y^2 +\varepsilon  Y(y)p_{x}^2}{X(x)-Y(y)}$.
 Then,  the geodesic flow of the 
 induced metric on the  quotient  space   $ \mathbb{R}^2/G$ (homeomorphic to the torus)   admits an integral quadratic in momenta.   We will call such metrics \emph{ globally--Liouville}.

\begin{figure}[ht!] 
  %\hspace{-2cm}
  {\includegraphics[width=0.7\textwidth]{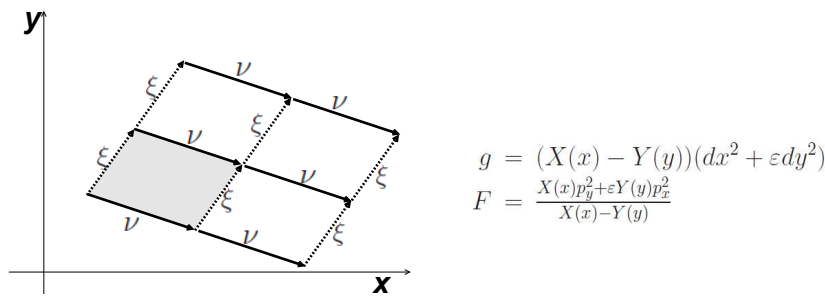}} 
  \caption{Vectors $\xi$ and $\nu$ and a fundamental region (gray parallelogram)  of the action of $G$ from Model Example 1. The torus $ \mathbb{R}^2/G$ can be identified with this   parallelogram with glued  opposite sides. Since the action of $G$ preserves $X(x)$ and $Y(y)$, the metric $g$ 
   induces a metric on  $ \mathbb{R}^2/G$, and the integral $F$  induces  an integral quadratic in momenta  }\label{liouville}
\end{figure}

\begin{Th} \label{main1}   Suppose a  metric $g$        on the two-torus  
${T}^2$ admits an  integral $F$  quadratic in momenta.  Assume  the integral is not  a linear combination of 
 the square of an integral linear in momenta and the Hamiltonian. Then, $(T^2, g)$ is globally Liouville,  i.e., there exist $X$, $Y$, $\xi$, $\nu$ satisfying the conditions in the Model Example 1 above   and a diffeomorphism $\phi:T^2 \to  \mathbb{R}^2/G$ that takes $g$ to the globally-Liouville metric $(X(x)- Y(y))(dx^2 +\varepsilon  dy^2)$ on  $\mathbb{R}^2/G$ and the integral $F$ to the integral $\pm  \left(\tfrac{X(x)p_y^2 +\varepsilon  Y(y)p_{x}^2}{X(x)-Y(y)}\right)$.  
\end{Th} 

In the Riemannian case, Theorem \ref{main1} follows from \cite{1,Kio}, see also \cite{BMF,BF}. 
We see that the answer in the pseudo-Riemannian case is essentially the same ( = no new phenomena appear) as the answer in the Riemannian case.  This similarity with the Riemannian case was unexpected: indeed, by Theorem \ref{main}, in the pseudo-Riemannian case (different from the Riemannian case) there are three  different types of metrics admitting quadratic integrals.  Moreover,  the examples from papers \cite{chanu,rastelli,smirnov,PR} show that, locally,  the pair (metric,integral) can change the type, i.e.,  the pair (metric,integral) can be, for example, as in Liouville case from one side of a line, and as in Complex-Liouville case from another side of the line. 
  But it appears that only one  type, namely the Liouville, 
 can exist on closed manifolds.

 Moreover, as we show in Example \ref{Ex2}, if the integral is the square of an integral linear in momenta, then the Jordan-block case is possible (even if the surface  is closed). Moreover, the pair 
 (metric,integral)  can change the type: be of Jordan-block type in a neighborhood of one point, and of Liouville type in a neighborhood of another point. Moreover, one can modify Example \ref{Ex2}(c) such that  the set of the points such that the pair (metric,integral)  changes the type is the direct product of the Cantor set and a circle.

 \subsection{Metrics on the Klein bottle  whose geodesic flows admit   integrals  quadratic in momenta }

  The scheme of the description is the same as for the torus: in Model Example 2 we describe a big family of metrics on the Klein bottle whose geodesic flows admit  integrals  quadratic in momenta. Theorem \ref{main3} claims that every metric such that its  geodesic flow admits an integral quadratic in momenta and such that the geodesic flow of the lift of the metric to the oriented cover admits no integral linear in momenta is as in Model Example 2. 
  \\[.1cm]

{\bf Model Example 2.} We consider $\mathbb{R}^2 $ with the standard coordinates $(x,y)$,  constants $c\ne 0, \ d\ne 0$,     two  vectors $\xi=(c,0)$ and $\nu=(0, d) $, and two nonconstant  functions $X$ and $Y$ of one variable (it is convenient to think that the variable of $X$ is $x$ and the variable of $Y$ is $y$)
 such that  \begin{itemize} 
 \item[(a)] for all  $(x,y)\in \mathbb{R}^2$ we have $X(x)\ne Y(y)$, and 
   \item[(b)] for every $(x,y)\in \mathbb{R}^2$, $ X(x+c)= X(x)$ and $ Y(y+d)= Y(-y)= Y(y)$. \end{itemize}
 Next, consider the metrics $(X(x)- Y(y))(dx^2 +\varepsilon  dy^2)$  on  $\mathbb{R}^2$ and 
 the action of the group   $G$
 generated by  the transformations $(x,y)\mapsto (x+ c, -y)$ and  $(x,y)\mapsto (x, y+ d)$.  The action is free, discrete and  preserves the metric and the quadratic integral  $ \tfrac{X(x)p_y^2 +\varepsilon  Y(y)p_{x}^2}{X(x)-Y(y)}$.
 Then,  the geodesic flow of the 
 induced metric on the  quotient  space   $ \mathbb{R}^2/G$ (homeomorphic to the Klein bottle)   admits an integral quadratic in momenta.   We will call such metrics \emph{ globally-(Klein)-Liouville}.   
 
 \begin{figure}[ht!] 
  %\hspace{-2cm}
  {\includegraphics[width=0.7\textwidth]{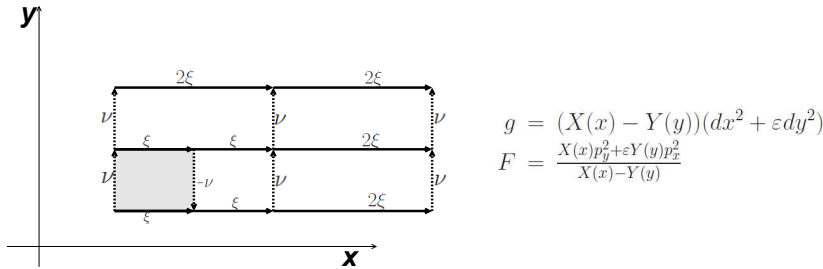}} 
  \caption{Vectors $\xi$ and $\nu$ and a fundamental region (gray rectangle)  of the action of $G$ from Model Example 2. The Klein bottle  $ \mathbb{R}^2/G$ can be identified with this   rectangle  with glued  opposite sides: the horizontal sides are glued with preserving the orientation, and the vertical sides are glued with inverting the orientation. The action of $G$ preserves  the metric $g$ 
    and the integral $F$; hence, the geodesic flow of the induced metric on  $ \mathbb{R}^2/G$ admits an  integral quadratic in momenta}\label{kleinliouville}
\end{figure}

 \begin{Th} \label{main3}   Suppose a  metric $g$        on the Klein bottle 
${K}^2$ admits an  integral $F$  quadratic in momenta.  Assume  the lift of the integral to the oriented cover is not   a linear combination of  the lift of the Hamiltonian and  
 the square of a function   linear in momenta. Then, $(K^2, g, F)$ is globally-(Klein)-Liouville,  i.e., there exist $X$, $Y$, $c,d$ satisfying the conditions in the Model Example 2 above   and a diffeomorphism $\phi:K^2 \to  \mathbb{R}^2/G$ that takes $g$ to the globally-(Klein)-Liouville metric $(X(x)- Y(y))(dx^2 +\varepsilon  dy^2)$ on  $\mathbb{R}^2/G$
 and  $F$ to the integral $\pm  \left(\tfrac{X(x)p_y^2 +\varepsilon  Y(y)p_{x}^2}{X(x)-Y(y)}\right)$.  
\end{Th} 
 
In the Riemannian case, Theorem \ref{main3} was proved  in \cite[Theorem 3]{klein}. We see that the answer in the pseudo-Riemannian case is essentially the same  as the answer in the Riemannian case (similar to the torus).

The following example explains why we require that the LIFT of the integral (to the oriented cover) is not   a linear combination of  the lift of the Hamiltonian  and    the square of a function   linear in momenta: 
 
 \begin{Ex} As in the Main Example 2,  we  consider $\mathbb{R}^2 $ with the standard coordinates $(x,y)$,  constants $c\ne 0, \ d\ne 0$,     two  vectors $\xi=(c,0)$ and $\nu=(0, d)$, the    function $X$  of the  variable $x$ such that $ X(x+c)= X(x)$.  Different from  the Main Example 2, by $Y$ we denote a CONSTANT 
 such that $X(x)\ne Y$ for all $x\in \mathbb{R}$.

 Under this assumptions, the metric  $(X(x)- Y)(dx^2 +\varepsilon  dy^2)$ and the integral $ \tfrac{X(x)p_y^2 +\varepsilon  Yp_{x}^2}{X(x)-Y}$ induce a  metric on the $K^2:= \mathbb{R}^2/G$, where $G$ is the group generated by the mappings $(x,y)\mapsto (x+ c, -y)$ and  $(x,y)\mapsto (x, y+ d)$, and an integral quadratic in momenta for the geodesic flow of this metric.

The lift of the integral to the oriented   cover $T^2:= \mathbb{R}^2/G'$, where $G':= \{ 2k\cdot \xi + m\cdot \nu \mid k, m\in \mathbb{Z}\}$,  is a linear combination of the Hamiltonian $\tfrac{1}{2} \frac{p_x^2 + \varepsilon p_y^2}{X(x)- Y}$ and the square of the (linear in momenta)  function $p_y$. Indeed,  $
F= p_y^2 + {2} \varepsilon Y  \cdot H. 
$
  
But, on ${K}^2$,   the integral in NOT a linear combination of the Hamiltonian and of the square of a    function  linear in momenta. The formal proof of this observation in the Riemannian case can be found in \cite[\S\S3,4]{klein}, the Riemannian proof can be  easily   generalized (using  Theorem \ref{main1}  of our paper) to the pseudo-Riemannian metrics. The main idea of the proof is that the function $p_y$ does not generate a function on the  Klein bottle, since the mapping $(x,y)\mapsto (x+ c, -y)$
 changes the sign of this function.
\end{Ex} 
 
 \subsection{Metrics on the torus  whose    geodesic flows admit  integrals linear in momenta} 
 
 In order to complete the description of the metrics of signature (+,--) whose geodesic flows admit nontrivial integrals quadratic in momenta,  we need to describe the metrics of signature (+,--) 
 on the torus such that their geodesic flows admit nontrivial integrals linear in momenta. 
 
 In the Riemannian case, metrics with geodesic flows admitting integrals linear in momenta can be considered as a partial case of the metrics whose geodesic flows admit integrals quadratic in momenta. Indeed, up to an isometry, any such metric is essentially  as in Model Examples 1,2 (see \cite{BMF,BF}), the only difference is that the function $ X$ is constant. 
 In particular, it implies that one can always slightly perturb a metric  whose geodesic flow admits an integral linear in momenta such that the geodesic flow of the result admits an integral quadratic in momenta, but admits no integral linear in momenta. 
 
 It appears that in the pseudo-Riemannian case the situation is different. 
 
  Below, we construct  a  family of metrics on the torus whose  geodesic flows admit integrals linear in momenta. In Examples \ref{Ex2}, \ref{Ex3}, we use the construction to show that in the pseudo-Riemannian case the following new (compared with the Riemannian case)  phenomena appear:
  
  \begin{itemize}   
  \item Example \ref{Ex2}(a) shows  that metric and  the integral can be as in the Jordan-block case. 
  \item Example \ref{Ex2}(c) shows that  the metric and the integral can be as in the Jordan-block case 
  in one neighborhood and   as in the Liouville case in another neighborhood.  
  \item  Example \ref{Ex3} shows the existence of a metric whose  geodesic flow   admits an integral linear in momenta, such that no small perturbation of this metric  admits an integral quadratic in momenta which is not a linear combination of the square of  an integral linear in momenta and the Hamiltonian.
 \end{itemize}

 {\bf Construction.} 
 We consider  $\mathbb{R}^2$ with the standard coordinates $x,y$ and  the standard orientation,  the vector fields $\xi:=  (1,0)$,  $ \eta:=(0,1)$, and 
a smooth  foliation on $\mathbb{R}^2$ invariant with respect to the flow of the vector field $\xi$ and with respect to the mapping $(x,y)\mapsto (x,y+1)$. 
With the help of these  data, we construct  a metric of signature $(+,-)$ on $\mathbb{R}^2$  such that $\xi$ is a Killing vector field for this metric.  

\begin{figure}[ht!] 
  %\hspace{-2cm}
  {\includegraphics[width=0.7\textwidth]{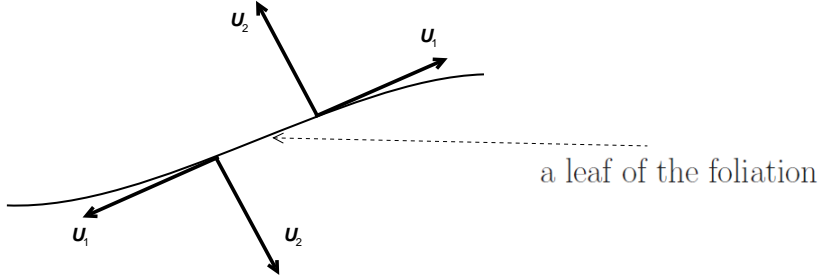}} 
  \caption{A leaf of the foliation and two possibilities for the vectors $U_1$, $U_2$}\label{u1u2}
\end{figure}

At every point $p$, we consider two vectors $U_1(p)$ and $U_2(p)$
 satisfying the following conditions: 
 \begin{itemize} 
\item $U_1$ at every point is tangent to the leaf of the foliation containing this point, 
\item $(U_1(p), U_2(p))$ is an orthonormal positive basis for the flat metric $\hat g =dx^2 + dy^2$, that is 
\begin{itemize}  \item[$\bullet$] $|U_1|_{\hat g}= |U_2|_{\hat g}=1$, $\hat g(U_1, U_2)=0$, 
\item[$\bullet$]   the orientation given by the basis coincides with the standard orientation, see Figure  \ref{u1u2}.\end{itemize}
 \end{itemize}
 
 Clearly, ar every point there exist precisely two possibilities for such vector fields $U_1$, $U_2$ (the second possibility is $(-U_1, -U_2)$).

 Now, consider the  metric $g$ such that in the  basis $(U_1,U_2)$ it has the matrix 
 $
 \begin{pmatrix} 
 0 & 1\\ 1& 0
 \end{pmatrix}.$ 
 The metric clearly 
 does not depend on the choice of vectors $U_1, U_2$ at every point, and is smooth. Since all objects we used to construct the metric are invariant with respect to the flow of $\xi$, the vector field $\xi$ is Killing for the metric. Then, the geodesic flow of the metric admits an integral $p_x$ linear in momenta.  
 Since all objects are invariant with respect to the lattice $G= \{ k\cdot \xi + m\cdot \eta\mid k,m\in \mathbb{Z}\}$, the metric induces  a metric on the torus $\mathbb{R}^2/G$ whose geodesic flow admits an integral linear in momenta.

   \begin{Rem} 
   By construction, the leaves of the foliation are light-line geodesics.
   \end{Rem}

\begin{figure}[ht!] 
  %\hspace{-2cm}
  {\includegraphics[width=1.1\textwidth]{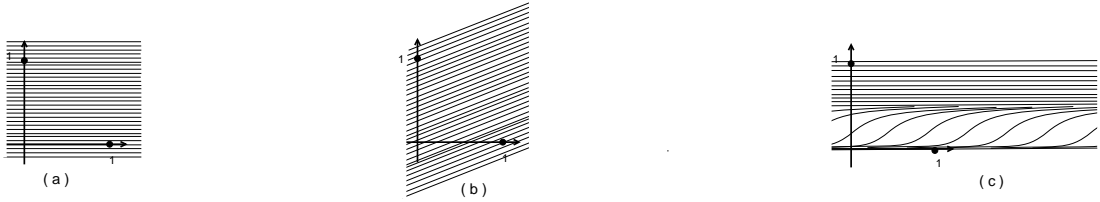}} 
  \caption{The foliations from Example \ref{Ex2}}\label{fig1}
\end{figure}

 \begin{Ex}\label{Ex2}   If the foliation is as on Figure  \ref{fig1}(a), the square of the integral is as in the Jordan block case. If the foliation is as on Figure  \ref{fig1}(b), the square of the integral is as in the Liouville case. If the foliation is as on Figure  \ref{fig1}(c), the square of the integral is as the Jordan block case in an annulus $\{(x,y)\in \mathbb{R}^2 \mid y-[y]> \tfrac{1}{2}\}$ and as in  Liouville case in the annulus  $\{(x,y)\in \mathbb{R}^2 \mid y-[y]<\tfrac{1}{2}\}$, where $[y]$ denotes the integer  part of $y$.   
\end{Ex}

\begin{center}
\begin{figure}[ht!] 
  %\hspace{-2cm}
  {\includegraphics[width=.5\textwidth]{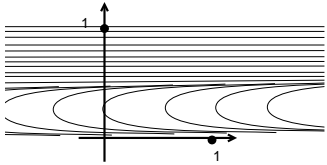}} 
  \caption{The foliation from Example \ref{Ex3}}\label{fig2}
\end{figure}
\end{center}

 \begin{Ex}\label{Ex3}   Let  the foliation is as on Figure  \ref{fig2}  (the restriction of the  foliation  to the annulus $\{(x,y)\mid x-[x]<\tfrac{1}{2}\}$  is the so-called Reeb component).     Then, the geodesic flow of 
  no  small perturbation of this metric admits an   integral quadratic in momenta that  is not a linear combination of the Hamiltonian and the square of an integral linear in momenta. Indeed,   the  Reeb component is stable with respect to small perturbations, and the light line geodesics of the 
  metrics from Model Example 1  are winding on the torus and form no Reeb component. 
\end{Ex}

Let us now describe all metrics on closed manifolds whose geodesic flows admit nontrivial integrals linear in momenta.  

{\bf Model Example 3.} 
We consider $\mathbb{R}^2 $ with the standard coordinates $(x,y)$,   the vectors $\xi:=(1,0)$ and $\nu:=(0, 1)$, and three  functions $K(y), L(y),M(y)$ of the variable $y$  periodic with period $1$  such that at every point $\det\begin{pmatrix}K & L \\ L&  M \end{pmatrix}=  KM-L^2<0$.  
 Next, consider the metric  $g= K(y)dx^2 +2L(y) dxdy + M(y)  dy^2$  on  $\mathbb{R}^2$,   and 
 the action of the lattice  $G:=\{k \cdot \xi  + m \cdot \nu \mid k, m\in \mathbb{Z}\}$ on $\mathbb{R}^2$. The action is free, discrete and   preserves the metric and the integral   $ p_{x}$ linear in momenta.
 Then,  the geodesic flow of the 
 induced metric on the  quotient  space   $ \mathbb{R}^2/G$ (homeomorphic to the torus)   admits an integral linear  in momenta.

\begin{Th} \label{main4}  Let   $g $ be a metric of signature $(+,-)$ on the torus $T^2$ such that it is not flat. 
If the geodesic  flow of $g$ admits an integral linear in momenta, then   the metric is 
 as in Model Example 3, i.e., there exist functions $K(y), M(y), L(y)$  periodic with period
 $1$ and a diffeomorphism $\phi:T^2\to \mathbb{R}^2/G$ that takes the metric $g$ to the metric  $ K(y)dx^2 +2L(y) dxdy + M(y)  dy^2$, and the integral to  $\const \cdot p_x$.  
\end{Th} 

In Theorem \ref{main4}, we assume that the metric $g$ is not flat. For  flat metrics, Theorem \ref{main4} is wrong, since the integral curves  of the Killing vector field corresponding to the linear integral  are non necessary closed curves for the flat metrics, but are closed curves in  Model Example  3. 
We need therefore to describe separately  flat metrics of  signature (+,--) on the torus.

By the { \it standard flat torus } we will consider $(\mathbb{R}^2/G, dxdy)$, where $(x,y)$ are    the standard coordinates on $\mathbb{R}^2$, and $G$ is a lattice generated by two  linearly independent vectors.

In \S \ref{flat} we will recall  why   every  torus $(T^2, g)$ such that the metric $g$  is flat and has signature (+,--) is isometric to a  standard one. 

\section{Applications} 
\subsection{Application I: Betrami problem on closed pseudo-Riemannian manifolds}  \label{beltrami} 
 Two metrics $g$ and $\bar g$ on one manifold are {\it geodesically equivalent,} if every  (unparametrized) geodesic of the first metric is a geodesic  of the second metrics. Investigation of geodesically equivalent metrics  is a classical topic  in differential geometry, see the surveys \cite{Aminova2,mikes}  or/and the introductions to \cite{threemanifolds,hyperbolic,diffgeo}. 
 
 In particular, in 1865 Beltrami  \cite{Beltrami} asked\footnote{ Italian original from \cite{Beltrami}: 
La seconda $\dots$  generalizzazione $\dots$ del nostro problema,   vale a dire:   riportare i punti di una superficie sopra un'altra superficie in modo  che alle linee geodetiche della prima corrispondano linee geodetiche della seconda.} {\it  to describe all pairs of  geodesically  equivalent Riemannian  metrics on  surfaces. } 
From the context it is clear that he considered this problem locally, in a neighbourhood of almost every  point, but the problem has sense, and is even more interesting globally.

Geodesically equivalent metrics and quadratic integrals are closely related: 
\begin{Th} \label{main2}  
Two metrics $g$ and $\bar g$ on $M^2$  are geodesically equivalent, if and only if  the following (quadratic in momenta) function  \begin{equation}\label{Integral}  F:TM^2\to \mathbb{R}, \
\
 F(x_1, x_2, p_1, p_2):= \left(\frac{\det(g)}{\det(\bar g)}\right)^{2/3}\cdot \sum_{i,j} \bar g^{ij}p_ip_j, 
\end{equation}
where we raised the indexes of $\bar g$ with the help of $g$, i.e., $\bar g^{ij} = g^{ki} \bar g_{km} g^{mj}$, 
 is an integral of the geodesic flow of  $g$. Moreover, $F=\const\cdot  H $ for  a certain $\const \in \mathbb{R}$ if and only if $g$ and $\bar g$ are proportional  with a  constant coefficient of proportionality. 
\end{Th}

  Theorem \ref{main2} above  was essentially known to  Darboux~\cite[\S\S600--608]{Darboux}; for recent proofs  see  \cite[Corollary 1]{pucacco}. See also the discussion  in  \cite[Section 2.4]{bryant}.

Combining  Theorems \ref{main1}, \ref{main3},  \ref{main4} with Theorem   \ref{main2}, we obtain a complete description of geodesically equivalent  pseudo-Riemannian metrics on closed  surfaces.

 \subsection{ Application II: every  quadratically-superintegrable metric  on a  closed surface  has constant curvature  } \label{superi}
 
 Recall that a metric on $M^2$ 
  is called quadratically-superintegrable, if the geodesic flow of the metric admits three linearly independent
   integrals quadratic  in momenta.  Quadratically-superintegrable metrics were first considered by Koenigs \cite{koenigs}.  Nowdays, investigation of quadratically-superintegrable metrics is  a hot topic in mathematical physics  due to various applications and  deep mathematical structures behind it, see e.g. \cite{kress}.

 For example, the standard  flat metric $dxdy$ on the 2-torus $\mathbb{R}^2/G$, where $G$ is a lattice generated by two linearly independent vectors, is quadratically-superintegrable. Indeed, the  Hamiltonian $H= 2p_xp_y$ and the quadratic in momenta functions $F_1:= p_x^2, \ F_2:= p_y^2$ are linearly independent integrals, and are invariant with respect to any lattice. 
 
 \begin{Cor} \label{super} 
 Let a metric $g$ on a closed surface  be  quadratically-superintegrable. Then, it has constant curvature.
 If in addition the metric has signature  {\rm (+,--)}, then   it is flat. 
 \end{Cor} 
 
 In the proof of Corollary \ref{super} we will need the following 
 
 \begin{Lemma} \label{-5} Let the metric $g$ of signature $(+, -)$ on the two-torus $T^2$ admit an integral quadratic in momenta that is not a linear combination of the Hamiltonian and of the square of an integral linear in momenta. Then, there exists a Riemannian  
 metric $\bar g$ geodesically equivalent to $g$. 
 \end{Lemma} 
 
 {\bf Proof.} By Theorem \ref{main1}, without loss of generality we can assume that  the metric $g$ and the integral $F$ are   as in Model Example 1. Without loss of generality we can think that $X(x)>Y(y)$ for all $(x,y)\in \mathbb{R}^2$. 
 
 Let us cook with the help of $H, F_1$ a Riemannian metric $\bar g$ geodesically   equivalent to $g$. 
 We put $X_{min}= \min_{x\in\mathbb{R} }X(x)$  and   $Y_{max}= \max_{y\in\mathbb{R} }Y(y)$. Clearly,  $X_{min}>Y_{max}$.  
We consider 
$$\bar  F:= H + \frac{1}{X_{min} + Y_{max}} F_1= 
\frac{\tfrac{1}{2}- \tfrac{Y}{X_{min} + Y_{max}}}{X-Y} p_x^2 +  \tfrac{\frac{X}{X_{min} + Y_{max}}-\tfrac{1}{2}}{X-Y} p_y^2.    
$$
Since $X > \tfrac{X_{min} +Y_{max }}{2} >Y$, the integral $\bar F$ is positively  defined (considered as a quadratic form on $T^*M^2$). Consider the metric $\bar g$ constructed by $\bar F$ with the help of Theorem \ref{main2}.   The metric is positively defined (i.e., is Riemannian), and is geodesically equivalent to $g$. Lemma \ref{-5} is proved.

 {\bf Proof of Corollary \ref{super}.}  The Riemannian version of Corollary \ref{super} is known (see \cite[Theorem 5.1]{Kio} and \cite[Lemma 3]{klein}, see also \cite[Theorem 6]{CMH}).
Then, without loss of generality we can assume that the metric has signature (+,--). 
 
  Let $H, F_1, F_2$ be the linearly independent  integrals quadratic in momenta. If both $F_1$ and   $F_2$ are linear combinations of the square of  integrals  linear in momenta and the Hamiltonian,   the metric admits two Killing vector fields implying that it has constant curvature.

 Assume now that there exists 
  an integral quadratic in momenta that is not a linear combination of the Hamiltonian and of the square
   of an integral linear in momenta.  By Lemma \ref{-5}, there exists a Riemannian 
    metric $\bar g$    geodesically equivalent to $g$. 
    The metric $\bar g$  
  is also quadratically-superintegrable. Indeed, as it was proved  in \cite[Lemma 1]{dim2} (see also \cite[\S 2.8]{bryant} and \cite[Lemma 3]{kruglikov}), every metric geodesically equivalent to a quadratically-superintegrable metric is also quadratically-superintegable.  Then, by the Riemannian version of Corollary \ref{super} (which is known, as we recalled above), the metric $\bar g$ has constant curvature. Then, by the Beltrami Theorem (see \cite{Beltrami,beltrami-small}), the metric $g$ also has constant curvature. The first part of Corollary \ref{super} is proved. 
 
 If the metric has signature (+,--), then the surface if the torus or the Klein bottle. 
 By  the Gauss-Bonnet Theorem,   a metric of constant curvature on the torus or on the Klein bottle 
  is flat. Corollary \ref{super} is proved.

 \subsection{Application III:  Proof of projective Obata conjecture for two-dimensional pseudo-Riemannian metrics} 
Let $(M^n,g)$ be a pseudo-Riemannian manifold of dimension $n\ge 2$. 
Recall that a \emph{projective transformation} of  $M^n$  is a diffeomorphism of the manifold that takes unparameterized geodesics to geodesics.

The goal of this paper is to prove the two-dimensional  pseudo-Riemannian version of the following

{\bf Projective Obata  conjecture.} { \it   Let a connected  Lie group $G$  act on a closed 
 connected   $(M^n, g)$ of dimension
 $n\ge 2$   by projective
transformations.  Then, it  acts by isometries, or  for some $c\in \mathbb{R}\setminus \{0\}$ the metric $c\cdot g$ is the Riemannian  metric of  constant positive sectional  curvature $+1$.}

\begin{Rem} The attribution  of conjecture to  Obata is in folklore (in the sense  we did not find a paper of Obata where he states  this conjecture). Certain papers, for example \cite{hasegawa,nagano,Yamauchi1},  refer to   this statement  as to  a classical conjecture.  If we replace ``closedness"  by 
``completeness", the  obtained conjecture is attributed in folklore to Lichnerowicz,  see also the discussion in  \cite{diffgeo}.  \end{Rem}

For Riemannian metrics, 
projective Obata conjecture was proved in \cite{obata,CMH,diffgeo}. Then,  in dimension two we may assume that the signature of the metric is $(+,  -)$, and that the manifold is covered by the torus $T^2$.  
  Thus, the two-dimensional version of the projective Obata conjecture follows from 
 
 \begin{Th} \label{obatath}  Let $(T^2, g)$ be the two-dimensional torus $T^2$ equipped with a metric $g$ of signature $(+, -)$. Assume 
  a connected  Lie group $G$  acts  on 
 $(T^2, g)$   by projective
transformations. 
 Then, $G$   acts by isometries.
\end{Th}

Note  that in the theory of geodesically equivalent metrics and projective transformations, dimension 2 is a special dimension: many methods that work in dimensions $n \ge 3$ do not work in dimension 2.  In particular, the proof of the projective Obata conjecture in the Riemannian case was separately done for dimension 2 in \cite{obata, CMH} and for dimensions greater than 2 in \cite{diffgeo}.  
Moreover, recently an essential progress was achived  in the proof of  the projective Obata conjecture in the pseudo-Riemannian case  in dimensions $n\ge 3$, see  \cite{kiosak2,mounoud}.   This progress allows us  to hope that it is possible to mimic (see \cite[\S 1.2]{kiosak2}) the Riemannian proof in  the pseudo-Riemannian situation (assuming the dimension is $n\ge 3$).  Thus, Theorem \ref{obatath}  closes an important partial case in the proof of projective Obata conjecture. 
 
 {\bf Proof of Theorem \ref{obatath}.}  Let $g$ be a pseudo-Riemannian metric  of signature $(+, -)$ of nonconstant curvature on $T^2$.  We denote by $\Proj_0(T^2, g)$  the connected component of the   group of projective transformations of $(T^2, g)$, and by $\Iso_0(T^2, g)$ the connected component of the group of isometries.   Clearly, $\Proj_0(T^2, g)\supseteq \Iso_0(T^2, g)$; our goal is to prove $\Proj_0(T^2, g)=\Iso_0(T^2, g)$. 
 
   We assume that  $\Proj_0(T^2, g)\ne \Iso_0(T^2, g)$.
     Then, there exists  a vector field  $v$  such that 
     it is  a projective vector field, but is not   Killing vector field.  (Recall that 
      a vector field $v$  is {\em projective}, if its local flow takes geodesics considered as unparameterized curved 
 to geodesics).  
    Then,  by \cite[Korollar 1]{obata},  \cite[Corollary 1]{CMH}, or \cite{Topalov}, 
  the quadratic in velocities function    
$$
 I:TM\to \mathbb{R}, \ \  \ I(\xi):=({\cal L}_vg)(\xi, \xi)-\tfrac{2}{3} \mathbb{\rm
trace}(g^{-1}{\cal L}_vg)\, g(\xi, \xi),  $$  where $\mathbb{\rm
trace}(g^{-1}{\cal L}_vg):= g^{ij} (\mathcal{L}_vg)_{ij}$ 
is a nontrivial (i.e., $\ne 0$)   integral  
for   the geodesic flow of  $g$.

Suppose first 
 $I$  is not a linear combination of the energy integral $g(\xi, \xi)$ and of the square of an integral linear in velocities.  Since  closed manifolds do not allow vector fields $v$ such that $\mathcal{L}_vg = \const\cdot  g $ for $\const \ne 0$, $I$ is  not proportional to the energy integral $g(\xi, \xi)$. 
 Then, by Lemma \ref{-5}, there exists a RIEMANNIAN metric $\bar g$ geodesically equivalent to $g$.

Every projective vector field for $g$ is also a projective vector field for $\bar g$ and vice versa, so that $\Proj_0(M, g) = \Proj_0(M, \bar g)$.  
 By  the (already proved)  Riemannian version of projective Obata conjecture we obtain that   $\Iso_0(M, \bar g) = \Proj_0(M, \bar g)$. Thus, $\Proj_0(M, g)= \Iso_0(M,\bar g).$   

  By \cite[Corollary 1]{beltrami-small}, see also \cite{knebelman},  the dimensions of the Lie group of isometries  of geodesically equivalent metrics coincide. Indeed,  for every Killing vector field $\bar K$ for $\bar g$ the vector field 
  $K^i:= \left(\tfrac{\det  g }{\det \bar g}\right)^{\tfrac{1}{n+1}}\bar g^{ik}g_{kj} \bar K^j $
is  a Killing vector field for $g$.  Then,  $\dim(\Iso_0(M, g)) = \dim(\Iso_0(M, \bar g))$ 
implying that $\Iso_0(M, g) =\Proj_0(M, g)$.   Hence, 
the assumption that $I$  is not a linear combination of the energy integral $g(\xi, \xi)$ and of the square of an integral linear in velocities  leads to a contradiction.    Thus, there exists a nontrivial integral  linear in velocities. Finally, there exists a nontrivial Killing vector field  that we denote by $K$.

Then, the group $\Proj_0$ is at least two-dimensional (because it algebra contains  $K$ and $v$).  
The structures of possible Lie groups of projective transformations was understood already by S. Lie \cite{Lie}. He proved that  the for a  $2-$dimensional 
 metric  of nonconstant curvature the Lie algebra of $\Proj_0$ is  the noncommutative 
  two dimensional algebra, or is $\eusl(3,\bbR)$. In both cases there exists a projective vector field 
  $u$ such that the linear span $span(u, K)$ is a two-dimensional noncommutative Lie algebra. 
  Then, without loss of generality we can assume that $[K, u]= u$ or $[K, u]= K$.

Now, by Theorem \ref{main4}, there exists a global coordinate system  $\bigl(x\in \ (\mathbb{R}, \ \textrm{mod} \ 1), y\in \ (\mathbb{R}, \ \textrm{mod} \ 1)\bigr)$  such that in this coordinate system 
$K= \alpha \cdot \tfrac{\partial }{\partial x}$, where $\alpha\ne 0$. Assume $u(x,y)=u_1(x,y) \tfrac{\partial }{\partial x} +u_2(x,y) \tfrac{\partial }{\partial y} $.  Without loss of generality we assume that $(u_1(0,0), u_2(0,0))\ne (0,0)$.  

Let $\phi_t$ be the  flow of $K$. Since $K= \alpha \cdot \tfrac{\partial }{\partial x}$, $\phi_t(x,y)=(x+ \alpha t, y)$. Let us calculate  the vector  $d\phi_t(u(0,0))$ for $t=1/\alpha$ by two methods (and obtain two different results which gives us a contradiction). 

First of all, since  $\phi_{1/\alpha}$ is the identity diffeomorohism,  $d\phi_t(u(0,0)) =u(0,0)$ for $t=1/\alpha$. 

The other method of calculating  $d\phi_t(u(0,0))$ is based on the commutative relation $[K, u]= u$ or $[K, u]= K$.

Let us first assume that $K,u$ satisfy  
$[K, u]= u$. In the coordinates, this  condition 
 reads  $ \alpha \tfrac{\partial }{\partial x}u_1 = u_1$ and  $ \alpha \tfrac{\partial }{\partial x}u_2 = u_2$ implying
$u_1(x,0) =  u_1(0,0) \cdot e^{x/\alpha }$ and $u_2(x,0) =  u_2(0,0) \cdot e^{x/\alpha }$.  Then,  $$d\phi_{1/\alpha}(u(0,0))=  u_1(0,0) \cdot e^{1/\alpha^2 }  \tfrac{\partial }{\partial x} +   u_2(0,0) \cdot e^{1/\alpha^2}  
 \tfrac{\partial }{\partial y} =  u(0,0)  \cdot e^{1/\alpha^2} . $$  Since $(u_1(0,0),u_2(0,0))\ne (0,0)$ we obtain  that $d\phi_{1/ \alpha}(u(0,0)) \ne u(0,0)$  which gives a contradiction.
  Thus, the commutative relation $[K, u]= u$  is not possible. 
 
 Let us now consider the second possible commutative relation $[K, u]= K$. 
 In coordinates  this   relation  reads 
 $ \alpha \tfrac{\partial }{\partial x}u_1 = \alpha$ and  $ \alpha \tfrac{\partial }{\partial x}u_2 = 0$ implying
$u_1(x,0) =  u_1(0,0) + {x}$. We again obtain that $d\phi_{1/\alpha}(u(0,0)) \ne   u(0,0)$,  which  gives   a contradiction.  Thus, the commutative relation $[K, u]= K$  is also not possible.  Finally, in all cases  the existence of a nontrivial projective vector field  on the torus 
 $T^2$ equipped with a metric of nonconstant curvature 
 leads to a contradiction. 
 
Let us now consider the remaining case: we assume that $g$ has constant curvature. 
By Gauss-Bonnet Theorem, a   metrics of constant curvature on $T^2$ is flat. Then, as we show in  \S \ref{flat},  $(T^2, g)$ is isometric 
to  the { \it standard flat torus } $(\mathbb{R}^2/L, dxdy)$, where $(x,y)$ are    the standard coordinates on $\mathbb{R}^2$, and $L$ is a lattice generated by two  linearly independent vectors.   In particular, all geodesics of the lift of the metric to $\mathbb{R}^2$ are  the standard straight lines.  Clearly, any projective
transformation of $(\mathbb{R}^2/L, dxdy)$ generates a bijection $\phi:\mathbb{R}^2 \to \mathbb{R}^2 $ that commute with the lattice $L$ and maps straight lines to  straight lines. It is easy to see that 
the connected component of the group of such bijections  consists of parallel translations, i.e., acts by isometries. Finally, $\Proj_0(\mathbb{R}^2/L, dxdy) = \Iso_0(\mathbb{R}^2/L, dxdy)$. Theorem \ref{obatath} is proved.

\weg{

{\bf Acknowledgement.}   The author  thanks
Deutsche Forschungsgemeinschaft
(Priority Program 1154 --- Global Differential Geometry and Research Training Group   1523 --- Quantum and Gravitational Fields)   and FSU Jena  for partial financial support, and D. Alekseevsky, O. Bauer, A. Bolsinov, G. Manno,  P. Mounoud,   G. Pucacco, and A. Zeghib  for useful discussions.  }

\section{Local theory and  the proof of Theorem~1 }  \label{puc}

\subsection{Admissible coordinate systems and 
Birkhoff-Kolokoltsov forms} \label{admissible}

Let $g$ be a pseudo-Riemannian metric   of signature (+,--) on connected  oriented  $M^2$.
Consider (and fix)  two vector fields $V_1, V_2$ on $M^2$  such that
\begin{itemize}
 \item[(A)]  $g(V_1, V_1) =g(V_2, V_2)=0$ and
 \item[(B)] $g(V_1, V_2)>0$, 
 \item[(C)] the basis $(V_1, V_2)$ is positive (i.e., induces the positive orientation). 
 \end{itemize}
Such vector fields always exist locally. Since locally there is precisely  two possibilities in choosing the directions of such vector fields, the  vector fields exist 
     on a finite (at most, double-) cover of $M^2$.

We will say that a local  coordinate system $(x,y)$
 is {\it admissible}, if  the vector fields  $\frac{\partial }{\partial x}$ and  $\frac{\partial }{\partial y} $ are proportional to $V_1, V_2$ with positive coefficient of proportionality:  $$\frac{\partial }{\partial x}=  \lambda_1(x,y) V_1(x,y), \  \  \  \frac{\partial }{\partial y}=  \lambda_2(x,y) V_2(x,y), \  \  \ \textrm{where $\lambda_i>0$}.$$
 Obviously,
\begin{itemize}
\item  admissible coordinates exist in a sufficiently small neighbourhood of every point, \item the metric $g$  in  admissible coordinates has the form
 \begin{equation}\label{metric}
 g =f(x,y)dxdy , \  \  \ \textrm{where $f>0$},  \end{equation}
     \item two admissible  coordinate systems  in
      one neighbourhood  are connected by \begin{equation} \label{coordinatechange} \begin{pmatrix} x_{new}\\
 y_{new}\end{pmatrix}
 = \begin{pmatrix} x_{new}(x_{old}) \\
 y_{new}(y_{old})\end{pmatrix} , \  \ \textrm{where  $\frac{dx_{ new}}{dx_{old}}>0$, $\frac{dy_{ new}}{dy_{old}}>0$}. \end{equation}
\end{itemize}

\begin{Rem} \label{admi} For further use let us note that smooth local  functions  $x,y$  form an  admissible coordinate system, if and only if $V_1(x)>0$, $V_2(y)>0$, and $V_1(y)= V_2(x)=0$ (where $V(h)$ denotes the derivative of the function $h$ in the direction of  the vector $V$). 
\end{Rem} 

\begin{Lemma}[\cite{pucacco}]  \label{BK}
Let $(x,y)$ be an admissible coordinate system for $g$.
Let  $F$ given by  \eqref{integral} be  an  integral for  $g$.
Then,
$$
B_1:= \frac{1}{\sqrt{|a(x,y)|}}dx, \;\; \left({\rm respectively },
B_2:= \frac{1}{\sqrt{|c(x,y)|}}dy \right)
$$
is a    1-form, which is defined at points such that $a\ne 0$ (respectively, $c\ne 0$).  Moreover, the coefficient
 $a$ (respectively, $c$) depends only on $x$ (respectively, $y$), which in particular implies that the forms $B_1$, $B_2$ are  closed.   \end{Lemma}

\begin{Rem}  The forms $B_1, B_2$ are not the direct analog of the ``Birkhoff" 2-form introduced by
Kolokoltsov in \cite{Kol}. In a certain sense, they are  real
analogs of the  two branches of the  square root of the Birkhoff  form.
\end{Rem}

\noindent{\bf Proof of Lemma~\ref{BK}.} The first part of the statement, namely that $$
\frac{1}{\sqrt{|a(x,y)|}}dx, \;\; \left({\rm respectively },
\frac{1}{\sqrt{|c(x,y)|}}dy \right)
$$
transforms as a $1$-form under admissible coordinate changes is evident:  indeed, after the coordinate change
\eqref{coordinatechange}, the  momenta transform as follows:
 $p_{x_{old}}= p_{x_{new}}\frac{d{x_{new}}}{d{x_{old}}}$, $p_{x_{old}}= p_{x_{new}}\frac{d{x_{new}}}{d{x_{old}}}$. Then, the integral $F$ in the new coordinates has
  the form
  $$ \underbrace{\left(\frac{d{x_{new}}}{d{x_{old}}}\right)^2{a}}_{a_{new} } {p_{x_{new}}^2}  + \underbrace{\frac{d{x_{new}}}{d{x_{old}}}\frac{d{y_{new}}} {d{y_{old}}}{b}}_{b_{new}}  {p_{x_{new}}} {p_{y_{new}}} + \underbrace{\left(\frac{d{y_{new}}}{d_{y_{old}}} \right)^2{c}}_{c_{new}} {p_{y_{new}}^2}.$$
  Then, the formal expression $\frac{1}{\sqrt{|a|}}dx_{old}
  $ ({\rm respectively}, $ \frac{1}{\sqrt{|c|}}dy_{old}$)   transforms into
  $$
  \frac{1}{\sqrt{|a|}} \frac{d{x_{old}}}{d{x_{new}}} dx_{new} \ \ \ \ \  \left(\textrm{respectively,  $ \frac{1}{\sqrt{|c|}}\frac{d{y_{old}}}{d{y_{new}}}dy_{new}$}\right),  $$
  which is precisely the transformation law of  1-forms.

Let us prove that  the coefficient
 $a$ (respectively, $c$) depends only on $x$ (respectively, $y$), which in particular implies that the forms $B_1$, $B_2$ are  closed.
 If $g$ is given by \eqref{metric}, its Hamiltonian is
 $$H=\frac{2p_xp_y}{f} \, ,$$
 and the condition $\{H, F\}=0$ reads \\
 \begin{eqnarray*}
0&=& \left\{\frac{2p_xp_y}{f}, ap_x^2+ bp_xp_y+ cp_y^2\right\}  \\
 &=& \frac{2}{f^2}\left(p_x^3(fa_y) + p_x^2 p_y (fa_x + fb_y + 2 f_x a + f_y b)+ p_yp_x^2 (fb_x + fc_y+ f_x b + 2 f_y c)+ p_y^3 (c_xf)\right) \, ,
 \end{eqnarray*}
 i.e., is equivalent to the following system of PDE:
 \begin{equation}\label{sys}
 \left\{\begin{array}{rcc} a_y&=&0 \, ,\\
 fa_x + fb_y + 2 f_x a + f_y b&=&0\, ,\\
 fb_x + fc_y+ f_x b + 2 f_yc &=&0\, ,\\
 c_x&=&0 \, .\end{array}
 \right.\end{equation}

 Thus, $a=a(x)$, $c=c(y)$ implying that 
  $B_1:= \frac{1}{\sqrt{|a|}} dx$ and $B_2:= \frac{1}{\sqrt{|c|}}dy$ are closed forms (assuming $a\ne 0$ and $c\ne 0$).  Lemma~\ref{BK} is proved.

\begin{Rem} \label{rem3} For further use let us formulate one more consequence of equations \eqref{sys}: if $a\equiv c \equiv 0$ in a neighbourhood of  a point, then $bf = \const$, implying     $F-    \tfrac{\const}{2} \cdot H=0$ in the neighborhood. If we consider   \eqref{sys}  as a system of PDE on the unknown functions $a,b,c$, we see that the system is linear and of finite type. Than, vanishing of the solution  corresponding to the integral $\hat F:= \left(F-  \tfrac{\const}{2} \cdot H\right)$ in the neighborhood implies vanishing of the solution on the whole connected manifold. Thus, if $a\equiv c \equiv 0$ in a neighborhood of  a point, then for a certain $\const \in \mathbb{R}$  we have $F\equiv   \const \cdot H$ on the whole manifold.
\end{Rem}

\begin{Rem} \label{rem4} For further use let us note that the set of the points where the form $B_1$ ($B_2$, resp.)  is not defined coincides with   the set of the points such that $a=0$ ($c=0$, resp.) and is invariant with respect to the (local) flow of the vector field $V_2$ ($V_1$, resp.)  
\end{Rem}

A local coordinate system $(x,y)$ will be called \emph{perfect}, if  it is admissible, and 
if in this coordinates system the coefficients  $a, c$ take values in the set $\{-1,0,1\}$ only.

\begin{Lemma} \label{adm} Let $F$ given by \eqref{integral} be an integral
for the geodesic flow of $ g =f(x,y) dxdy$ 
such that $F\ne \const \cdot H$ for all $\const \in \mathbb{R}$. 
Then,    almost every point $p$  has a neighborhood $U$ such that  precisely one of the following conditions is  fulfilled:  \begin{itemize}  \item[(i)]  $ac>0$ at all points of $U$, 
\item[(ii)]  $ac<0$ at all points of $U$, 
\item[(iii)(a)]  $a=0$ and  $c\ne 0$ at all points of $U$,   or 
\item[(iii)(b)]  $a \ne  0$   and $c=0$ at all points of $U$.
 \end{itemize} 

Moreover, there exists a perfect coordinate  system $\tilde x, \tilde y$   in a (possibly, smaller)   neighborhood $U'(p)\subseteq U(p)$  of $p$.
In the perfect coordinate system, the metric and the integral 
are given by 
$$
g= \tilde f(\tilde x,\tilde y)dxdy \ \ \textrm{and} \ \ F= \sign(a(x,y)) p_{\tilde x}^2 + \tilde b(\tilde x,\tilde y) p_{\tilde x}p_{\tilde y} + \sign(c( x,y)) p_{\tilde y}^2, 
$$
where $\sign(\tau) = \left\{\begin{array}{ccc} 1  &\textrm{if} & \tau >0 \\  -1  &\textrm{if} & \tau <0 \\ 0  &\textrm{if} & \tau =0.  \end{array}\right.$ 
\end{Lemma} 

\noindent {\bf Proof of Lemma  \ref{adm}.}  It is sufficient to prove the lemma  assuming that $M^2$ is a small neighborhood $W$. We consider and fix  admissible coordinates in this neighborhood. In this coordinates the coefficients $a,b,c$ of the integral \eqref{integral} are smooth functions.

We conisder the following subsets of $W:$ 
\begin{itemize}
           \item $W_{ac\ne 0}:= \{ q\in W \mid a(q)c(q) \ne 0\}$, 
           \item $W_{a\ne 0,c = 0}:= \{q\in W \mid a(q) \ne 0, c(q)= 0\}$, 
           \item $W_{a= 0,c \ne  0}:= \{q\in W \mid a(q) = 0, c(q)\ne  0\}$, 
\item $W_{a= 0,c = 0}:= \{q\in W \mid a(q) = 0, c(q)= 0\}$. 
\end{itemize}

The sets are clearly disjunkt, there union coincides with the whole $W$.  
We consider the set $W_{\textrm{perfect}}:= W_{ac\ne 0} \cup  int(W_{a=0,c \ne 0}) \cup  int(W_{a\ne 0,c = 0}), $ where ``$int$" denotes the set of inner points. The set $W_{\textrm{perfect}}$
 is  open, and is everywhere dense in  $W$. Indeed, it is open, since $W_{ac\ne 0} $,  $int(W_{a=0,c \ne 0})$, and $  int(W_{a\ne 0,c = 0}) $ are open.  It is everywhere dense,
  since it is  everywhere dense in the set   $W_{ac\ne 0} \cup  W_{a=0,c \ne 0}\cup W_{a\ne 0,c = 0}, $ and the remaining set $W_{a= 0,c = 0}$  is nowhere dense by  Remark \ref{rem3}. 
 
Now, by definition,  every point of $W_{\textrm{perfect}}$  has a neighborhood such that in this neighborhood one of the conditions (i)--(iii) is fulfilled.   The first statement of the proposition is proved. 

Let us now 
 prove the second statement.   
 Let  $p_0\in int(W_{a\ne 0,c = 0})$.  In a simply-connected 
neighborhood  $U(p_0)\subset W_{a\ne 0,c = 0}  $, we consider the function  
\begin{equation} \label{normalx}
 x_{new}(p) :=\int\limits_{p_0}^p
 B_1. \end{equation} Since the form $B_1$ is closed, and $U(p_0)$ is simply-connected, the function $x_{new} $ does not depend on  the choice of the  curve connecting the points $p_0,p$, 
 and is therefore well defined. The  differential of the   function  $x_{new}$ is precisely the 1-form $B_1$, and does  not vanish at $p_0$.  We have $V_1(x_{new})= B_1(V_1)>0$, 
 $V_2(x_{new})= B_1(V_2)= 0$. Since the coordinates $(x,y)$ are admissible,  $V_2(y)>0$ and 
 $V_1(y)=  0$.  Then, by Remark \ref{admi}, $(x_{new}, y)$ is a  local admissible 
 coordinate system in  a possibly smaller neighborhood $U'\subseteq U$ containing $p_0$. 
 
 \begin{Rem} Let us note that, in the admissible coordinates the formula \eqref{normalx} looks \begin{equation}\label{loc:normalx}  x_{new}(x{})=\int_{x_0}^{x_{}} \frac{1}{\sqrt{|a(t)|} }\, dt   
\end{equation}  implying that $x_{new}$ is independent of    $y$, i.e., $x_{new}=x_{new}(x)$. \end{Rem} 
 
  In this coordinate system,  
  the integral   $F$ is equal to $$ {\left(\frac{d{x_{new}}}{d{x}}\right)^2{a}} {p_{x_{new}}^2}  + {\frac{d{x_{new}}}{d{x_{old}}}{b}} {p_{x_{new}}} {p_{y}} = \frac{{a}}{(\sqrt{|a|})^2} {p_{x_{new}}^2}  + {\frac{b}{\sqrt{|a|}}}  {p_{x_{new}}} {p_{y}} = 
\sign(a) {p_{x_{new}}^2}  + b_{new}  {p_{x_{new}}} {p_{y}} .$$

  The cases $p_0\in int(W_{a= 0,c \ne  0})$, $p_0\in W_{a\ne 0,c \ne   0}$ are  similar: in the case $p_0\in int(W_{a =0,c \ne 0})$,  in the coordinate system $(x, y_{new})$ in a possibly smaller neighborhood of $p_0$, where 
  \begin{equation} \label{normaly}
 y_{new} :=\int\limits_{p_0}^pB_2,\end{equation}
the integral $F$ is given by $b_{new}p_{x}p_{y_{new}}  + \sign(c) p_{y_{new}}^2$. In the case  $p_0\in W_{a\ne 0,c \ne   0}$, 
 in the coordinate system $(x_{new}, y_{new})$, where $x_{new}$ is given by \eqref{normalx} and $y_{new}$ is given by \eqref{normaly}, the integral $F$  is given by  $\sign(a) p_{x_{new}}^2+ b_{new}p_{x_{new}}p_{y_{new}}  + \sign(c) p_{y_{new}}^2$. Lemma   \ref{adm}
 is proved.

\begin{Rem}  \label{pp} If $a=0$ ($c=0$, resp.), the coordinate 
transformation of the form  $(x_{new}(x),y)$ ($(x,y_{new}(y))$, resp.)  does not change the property of coordinates to be perfect.  
 If $ac\ne 0$, 
the perfect  coordinates are unique up to  transformation $(x,y)\mapsto (x+ \const_1, y+\const_2)$. In particular, if $ac\ne 0$, the vector fields $\tfrac{\partial }{\partial x}$ and  $\tfrac{\partial }{\partial y}$, where $x,y$ are local perfect  coordinates,
 do not depend on the choice of local perfect  coordinates, and therefore are well-defined  globally, at all points such that $ac\ne  0$ (provided that $V_1, V_2$ satisfying (A,B,C) are globally defined). \end{Rem}

\subsection{Proof of Theorem~2}

By Lemma  \ref{adm},  almost every point of $M^2$ has a neighborhood  such that 
 in perfect coordinates      the metrics and the integral are as in one of the following cases:

\begin{itemize} \item[] {\bf Case 1: $ac >0$:} The metric is $f(x,y) dx dy $, the  integral is $\pm (p_x^2 + b(x,y) p_xp_y + p_y^2)$. 

 \item[] {\bf Case 2: $ac<0$:}   The metric is $f(x,y) dx dy $, the  integral is $\pm(p_x^2 + b(x,y) p_xp_y - p_y^2)$. 

 \item[] {\bf Case 3a: $c\equiv 0$: } The metric is $f(x,y) dx dy $, the  integral is $\pm(p_x^2 + b(x,y) p_xp_y)$.  
  
   \item[] {\bf Case 3b: $a\equiv 0$: } The metric is $f(x,y) dx dy $, the  integral is $\pm( b(x,y) p_xp_y+ p_y^2)$.   \end{itemize} 
  We will carefully consider all four cases. 

  \subsubsection{Case 1} \label{case1} 
  
  \begin{Prop} \label{c1} Let the  geodesic flow of a metric $g=f(x,y)dxdy$  admits an integral \eqref{integral}. Assume $ac>0$  at the point $p$. 
   Then, in the coordinates $(u,v)= (\tfrac{x +y}{2}, \tfrac{x-y}{2})$, where $(x,y)$ are perfect  coordinates in  a neighborhood of $p$,  
   \begin{equation}\label{answer:case1}
 g=(U(u)-V(v))(du^2 -dv^2) \  \textrm{and} \ \ 
 F= \pm \left(\frac{p_v^2 U(u) - p_u^2V(v)}{U(u)-V(v)}  \right)\, ,
 \end{equation}
 where $U,V$ are certain functions of one variable.
  \end{Prop} 
  
\noindent{\bf    Proof. } Without loss of generality $a$ and $c$ are positive in a neighborhood of $p$.  
 Then, by  Lemma \ref{adm}, in  perfect  coordinates  in a neighborhood of   $p$
  the metric and the integral are   
  $g=f(x,y) dx dy $, $F=p_x^2 + b(x,y) p_xp_y + p_y^2$. Then, the system (\ref{sys}) has the following simple form:
$$
 \left\{\begin{array}{rcc}  (fb)_y+ 2 f_x  &=&0 \,,\\
 (fb)_x  + 2 f_y&=&0 \, ,  \end{array}
 \right. \ \textrm{which is equivalent to} \ 
 \left\{\begin{array}{rcc}  (fb+ 2 f)_x  + (fb  + 2 f)_y&=&0 \,,\\
   (fb  -2 f)_x  -(fb- 2 f)_y&=&0 \,. \end{array}
 \right.$$
 After the (non-admissible) change of coordinates $u = \tfrac{x+y}{2}$, $v= \tfrac{x-y}{2}$,  the system has  the form
$$
 \left\{\begin{array}{rcc}  (fb+ 2 f)_u&=&0 \,,\\
   (fb  -2 f)_v&=&0 \,, \end{array}
 \right.\ \textrm{which is equivalent to} \ \left\{\begin{array}{rcc}   fb+ 2f & =&  -4{V(v)}\, ,\\  fb-2f& =&-4{U(u)}\end{array} \right.
 $$ for certain functions 
 $U(u) $  and $V(v)$. 
 Thus,
 $$
 f= {U(u)-V(v)} \, , \;\; b=-2 \frac{U(u)+V(v)}{U(u)-V(v)}\, .
 $$

 Let us now calculate the metric in the integral in the new coordinates:  substituting 
 $dx= du + dv, dy= du- dv$ in the formula $g = f(x,y) dxdy = {(U(u)-V(v))} dxdy$, we obtain that  in the new coordinates the metric is $(U(u)-V(v))(du^2 - dv^2)$. Substituting $p_x =\left(\tfrac{\partial u}{\partial x}p_u  + \tfrac{\partial v}{\partial x}p_v \right) = \tfrac{1}{2} (p_u + p_v)  $ and   
 $p_y =\left(\tfrac{\partial u}{\partial y}p_u  + \tfrac{\partial v}{\partial y}p_v \right)=\tfrac{1}{2} (p_u - p_v) $ in the formula $F= p_x^2  + b p_xp_y +p_y^2 = p_x^2  -2 \frac{U(u)+V(v)}{U(u)-V(v)}p_xp_y +p_y^2$, we obtain that  in the new coordinates $(u,v)$ 
 $$  
  F=  \tfrac{1}{2}\left(  p_u^2 + p_v^2  - \frac{U(u)+V(v)}{U(u)-V(v)}(p_u^2 - p_v^2)  \right)=
     \frac{U(u) p_v^2 - V(v) p_u^2         }{U(u)-V(v)}.  
 $$
We see that,  in the  new coordinates,  the metric and the integral   are as in \eqref{answer:case1}. Proposition \ref{c1} is proved.

  \subsubsection{ Case 2  }

  \begin{Prop} 
 \label{c2} Let the  geodesic flow of a metric $g=f(x,y)dxdy$  admits an integral \eqref{integral}. Assume $ac<0$  at the point $p$.  Then, in  perfect coordinates in a neighborhood of $p$,  
   \begin{equation}\label{answer:case2}
g=  \Im(h)dxdy \ \ \ \textrm{and} \ \ \ F=\pm \left(  p_x^2 - p_y^2   + 2\frac{\Re(h)}{\Im(h)}p_xp_y\right),
 \end{equation}
where $\Re(h)$ and $\Im(h)$ are the real and the imaginary  parts of a holomorphic function $h$ of the variable $z = x+ i\cdot y$.  
  \end{Prop} 
\noindent{\bf    Proof. }
Without loss of generality $a(p)>0, $ $c(p)<0$. By Lemma   \ref{adm}, in  perfect coordiantes 
  the metric and the integral are   
  $g=f(x,y) dx dy $, $F=p_x^2 + b(x,y) p_xp_y - p_y^2$. Then, the system (\ref{sys}) has the following simple form:
  
 \begin{equation}\label{sys:case2}
 \left\{\begin{array}{rcc}  (fb)_y+ 2 f_x  &=&0 \,,\\
 (fb)_x  -2 f_y&=&0 \,. \end{array}
 \right.\end{equation}
 We see that these equations are the Cauchy-Riemann  conditions for the complex-valued
 function $fb+ 2i f$. Thus, for an appropriate holomorphic
 function $h= h(x+ iy)$ we have $fb=\tfrac{1}{2}\Re(h)$, $f =\Im(h)$.
 Finally,  the metric and the integral  have the form \eqref{answer:case2}.   Proposition  \ref{c2} is proved. 
 
 \subsubsection{Case 3  } \label{3t} 

In this case we prove  two propositions: the first one is more general, and is the final step in the proof of Theorem \ref{main}. The second one   requires additional assumptions, and will be used in the proof of Theorem  \ref{main1}.  

 \begin{Prop} 
 \label{c3} Let the  geodesic flow of a metric $g=f(x,y)dxdy$
  admits an integral \eqref{integral}.  Then, the following two statements are true: 
 
 \begin{itemize}   
  \item[(a)] If  $a(p)\ne 0$, and $c(q)=0$    at every point $q$  of a small neighborhood  of  $p$,   in  perfect  coordinates in a (possibly, smaller) neighborhood of $p$,  
   \begin{equation}\label{answer:case3a}
g=  \left( \widehat{ Y}(y)+\frac{x}{2} Y'(y)\right)dxdy \ \ \textrm{and} \  \ F= \pm\left(p_x^2 - \frac{Y(y)}{\widehat{ Y}(y)+\frac{x}{2} Y'(y) }p_xp_y \right)\, ,
 \end{equation} where $Y$ and $\widehat Y$ are functions of one variable. 
 \item[(b)] If  $c(p)\ne 0$, and $a(q)=0$   at every point $q$  of a small neighborhood  of  $p$,   in  perfect  coordinates in a (possibly, smaller) neighborhood of $p$,  
   \begin{equation}\label{answer:case3b}
g=  \left( \widehat{ X}(x)+\frac{y}{2} X'(x)\right)dxdy \ \ \textrm{and} \  \  F=\pm \left( p_y^2 - \frac{X(x)}{\widehat{ X}(x)+\frac{y}{2} X'(x) }p_xp_y \right) \, ,
 \end{equation}  where $X$ and $\widehat X$ are functions of one variable.  \end{itemize}
\end{Prop}
\noindent{\bf    Proof. }
The cases (a) and (b) are clearly analogous;  
without loss of generality we can assume  $a(p)>0, $ $c\equiv 0$. 
 By Lemma \ref{adm}, in  perfect  coordinates 
  the metric and the integral are   
  $g=f(x,y) dx dy $, $ F=p_x^2 + b(x,y) p_xp_y $. 
  Then, the equation (\ref{sys}) has the following simple form:
 \begin{equation}\label{sys:case3}
 \left\{\begin{array}{rcc}  (fb)_y+ 2 f_x  &=&0 \, ,\\
 (fb)_x  &=&0 \, .\end{array}
 \right.\end{equation}
 This system can be solved. Indeed, the second equation implies $fb= -Y(y)$. Substituting this in the first equation  we obtain
 $Y'(y)= 2f_x$ implying
 $$f= \frac{x}{2} Y'(y)+ \widehat{ Y}(y) \textrm{ \ \ and}  \ \ \ b= - \frac{Y(y)}{\frac{x}{2} Y'(y)+ \widehat{ Y}(y)}\, .
 $$
 Finally, the metric and the integral are as in \eqref{answer:case3a}. Proposition \ref{c3}(a) is proved.  The proof of Proposition \ref{c3}(b) is essentially  the same.

{\bf Proof of Theorem 1. } 
  Theorem~\ref{main} follows directly from Lemma \ref{adm} and  Propositions \ref{c1}, \ref{c2}, \ref{c3}. Indeed,  
  by  Lemma \ref{adm}, almost every point has a neighborhood such that in this neighborhood the assumptions of one of Propositions \ref{c1}, \ref{c2}, \ref{c3} are fulfilled. Then, by Propositions \ref{c1}, \ref{c2}, \ref{c3} the metric and the integral are as in the table in  
 Theorem \ref{main}.
  
  We will also need a slightly less general version of normal  form of  metrics   satisfying the assumption of Case 3.   
  
  Let us observe  that the function 
  $Y$ from \eqref{answer:case3a}, or the function $X$  from \eqref{answer:case3b},  can be given in  invariant terms  (i.e.,   they  does not depend on the choice of a perfect coordinate system,   and 
  can be smoothly prolonged to the whole manifold). Indeed,  consider the symmetric 
  $(2,0)-$tensor $\tilde F^{ij}$ such that $F= \sum_{i,j}\tilde F^{ij}p_ip_j$ 
  (if $F $ is given by \eqref{integral}, the matrix of $\tilde F$ is $\begin{pmatrix} a & b/2 \\ b/2 & c\end{pmatrix}$).  Transvecting    $\tilde F^{ij}$ with $g_{ij}$ we obtain the globally 
  defined smooth function    $L:=  \trace(\tilde F_j^i):= \sum_{i,j}  \tilde F^{ij}g_{ij}$.  Under assumptions of Case 3a, in  the perfect   coordinates, the function $L$ is given by   
  \begin{equation}\label{L}
  L= \sum_{i,j}  \tilde F^{ij} g_{ij} =
   \trace\left(    \begin{pmatrix} a & b/2 \\ b/2 & c\end{pmatrix}  \begin{pmatrix} 0& f/2 \\ f/2 & 0\end{pmatrix}     \right) = \trace\left(    \begin{pmatrix} -Y/4& \ast \\ 0 & -Y/4\end{pmatrix}     \right)=-Y/2.
  \end{equation}

 \weg{ Now, let us observe that if  (under assumptions  of Case 3a) $Y\equiv \const$ in a neighborhood $U$, then in this neighborhood the integral $F$  is a linear combination of a square of an integral linear in momenta and the Hamiltonian. Indeed,  in this case  $Y'=0$, and the formulas \eqref{answer:case3a} look 
   $$g=   \widehat{ Y}(y)dxdy \ \ \textrm{and} \  \ F= \pm\left(p_x^2 - \frac{\const }{\widehat{ Y}(y)}  p_xp_y \right). 
$$  
  Since the components of the metric do not depend on $x$, $p_x$ is an integral (linear in momenta)
   of the geodesic flow of $g$, and 
  $F$  is a linear combination of    $p_x^2$ and $H$.

  Since in Theorem  \ref{main1} we require that the restriction of the integral to no neighborhood is a linear combination of a square of an integral linear in momenta and the Hamiltonian,  we can 
  take $L$ as a local coordinate in the neighborhood of almost every point $p$ such that $c\equiv 0$ in $U(p)$, or $a\equiv 0$ in $U(p)$.   }

 \begin{Prop} 
 \label{c4} Let the  geodesic flow of a metric $g=f(x,y)dxdy$
  admits an integral \eqref{integral}.  Then, the following two statements are true: 
 
 \begin{itemize}   
  \item[(a)] Suppose   $a(p)\ne 0$, and $c(q)=0$    at every point $q$  of a small neighborhood  of  $p$. 
   Assume $dL_{|p}\ne 0$, where $L$ is given by  \eqref{L}. Then,   in a (possibly, smaller) neighborhood of $p$,      in   perfect  coordinates  $(x,y) $  such that $y(q)=  -{2}L(q)$ for all  $q$, the metric and the integral are given by   
   \begin{equation}\label{answer:case3abis}
g=  \left( { Y}(y)+\frac{x}{2} \right)dxdy \ \ \textrm{and} \  \ F= \pm\left(p_x^2 - \frac{y}{{ Y}(y)+\frac{x}{2} }p_xp_y \right)\, , 
 \end{equation}
 where $Y$ is a  function of one variable. 
 \item[(b)] Suppose    $c(p)\ne 0$, and $a(q)=0$   at every point $q$  of a small neighborhood  of  $p$. 
  Assume $dL_{|p}\ne 0$, where $L$ is given by  \eqref{L}. Then, in a (possibly, smaller) neighborhood of $p$, 
    in  perfect  coordinates  $(x,y) $  such that $x(q)=  -{2}L(q)$ for all  $q$, the metric and the integral are given by   
   \begin{equation}\label{answer:case3bbis}
g=  \left( { X}(x)+\frac{y}{2} \right)dxdy \ \ \textrm{and} \  \ F= \pm\left(p_y^2 - \frac{x}{{ X}(x)+\frac{y}{2} }p_xp_y \right)\, , 
 \end{equation} 
 where $X$ is a  function of one variable.  \end{itemize}
\end{Prop}

  \noindent{\bf    Proof. } The cases (a) and (b) are clearly analogous;  
without loss of generality we can assume  $a(p)>0, $ $c\equiv 0$. In the perfect coordinates such that 
$y= -{2}L$,  we have  $g=f(x,y)dxdy$ and $F= p_x^2 -\tfrac{y}{f}  p_x p_{y} $. Then, the system    
 \eqref{sys} is equivalent to the  equation
 $
   2 f_x  =1.$
  Thus, $f=Y(y) +\tfrac{x}{2}$. Proposition \ref{c4}(a) is proved. The proof of Proposition \ref{c4}(b) is similar.

  \section{ Global theory and the main step in the  proof of Theorem~2}   \label{4} 
   
   \subsection{Notation, conventions, and the plan of the proof} \label{notation}
   
   Within the whole  section we assume that \begin{itemize} \item 
   the surface is the torus $T^2$,  
    \item $g$ is a pseudo-Riemannian metric of signature (+,--)  on $T^2$. 
   \item The vector fields $V_1$, $V_2$ satisfying conditions (A,B,C)  from \S \ref{admissible}  are globally defined (the case when it is not possible  will be considered in \S \ref{laststep}).  
    \item $F$  is a nontrivial  integral of the geodesic flow of
    $g$.   We will reserve  notation $x,y$ for  admissible coordinates, or for  perfect coordinates, and will denote the coefficients of the integral as in \eqref{integral}.   As in \S \ref{admissible}, we will denote by $B_1$, $B_2$  the $1-$forms 
   $\tfrac{1}{\sqrt{|a|}} dx$  and  $\tfrac{1}{\sqrt{|c|}} dy$.   \end{itemize} 
   
  As in \S \ref{3t}, we denote by $\tilde F^{ij}$ the 
   symmetric $(2,0)-$tensor corresponding to the integral $F$, and by $\tilde F_j^i$ the $(1,1)-$tensor 
   $\tilde F_j^i:= \sum_k \tilde F^{ik}g_{k j }.$ 
  
  We will proceed according to the following plan: 
  
  \begin{enumerate} 
  \item In \S \ref{inco} we show that there exists no point such that $ac<0$. This will imply  that  $\tilde F_j^i$ has real eigenvalues at every point of $T^2$.
   
  \item  By Remark  \ref{invar},  $\tilde F_j^i$  has only one    eigenvalue 
  (of algebraic multiplicity 2)  at the points such that $B_1$ or $B_2$  is not defined. In \S \ref{injo}, we show that 
  this eigenvalue is 
    constant on each connected component of the set such that $B_1$ or $B_2$ is not defined.  
  \item In \S  \ref{liou} we show that the existence a point such that $B_1$ or $B_2$ is not defined implies that 
  one of the eigenvalues of $\tilde F_j^i$ is constant  on the whole manifold. 
   
   \item In \S \ref{eigenvalue}, we show that if one of the eigenvalues of $\tilde F_j^i$ is constant, the quadratic integral $F$, or the  lift of the quadratic integral  to the appropriate double  cover is a linear combination  the square of a function  linear in momenta and the Hamiltonian. Later, in Corollary \ref{cover},
    we show that if the lift of the quadratic integral  to a double  cover is a linear combination of the  lift of the Hamiltonian  and  the square of a integral  linear in momenta,
    then the integral  is  a linear combination  of the Hamiltonian and the square of an integral linear in momenta. 
   
  \item In \S\ref{endof} we show that if at every point $B_1$ and  $B_2$  are defined, then  the torus, the 
   metric $g$,  and the integral $F$ are as in the Model Example 1.  
    \end{enumerate} 
   
   These will prove Theorem \ref{main1} under the additional assumption that the vector fields $V_1, V_2 $ exist on $T^2$. The case when this vector fields do not exist on $T^2$ will be considered later, in \S \ref{laststep}: we will prove that this  case can not happen (if there exists an integral quadratic in momenta that is not a linear combination of the Hamiltonian and the square of an integral linear in momenta). 
   
\subsection{At every point, the eigenvalues of $\tilde F_j^i$ are real }  \label{inco}

\begin{Lemma} \label{impcomplex} 
There is no point $p \in T^2$ such that at this point $ac<0$. 
\end{Lemma}

\noindent{\bf    Proof. }
Suppose at $p\in T^2$ we have $ac<0$. Let $W_0$ be the connected 
 component of the set 
$$W:=  \{ q \in T^2 \mid \textrm{$B_1$ and $B_2$ are defined} \}$$ 
containing the point $p$.  At every $q\in W_0$ we have $ac<0$. 
We consider the function $K:W_0\to \mathbb{R}$, $K= \frac{1}{g^*(B_1, B_2)}$, where $g^*$ is the scalar product on $T^*T^2$ induced by $g$. 

In any perfect  coordinates $(x,y)$  we have  $B_1= dx$,  $B_2=dy$, and $g = \Im(h) dx dy $ by Proposition \ref{c2}. Then,  $K= \Im(h)$ for a holomorphic   function $h$ implying it is harmonic function. When we approach the  boundary $\overline W_0\setminus W_0$, the function $K$ converges to $0$.  Indeed, in the admissible coordinates near a boundary point 
the function $K$ is $ f \sqrt{|ac|}$, and $ac\stackrel{\textrm{converges}}{\longrightarrow}  0$ (because at least one of coefficients $a,c$ is zero at the points of boundary).

Finally, by the  maximum principle (for  harmonic functions), the function $h$ is identically zero, which clearly contradicts the assumptions.  Lemma \ref{impcomplex}  is proved.

\begin{Cor} \label{jjj}
At every point of $T^2$, the eigenvalues of $\tilde F_j^i$ are real.
\end{Cor} 
{\bf Proof.} The eigenvalues are the roots of the characteristic polynomial
\begin{equation} \label{2-1}\chi(t)= \det(\tilde F_j^i- t\cdot \delta_j^i)= \det\left(\begin{pmatrix} fb/4 & af/2 \\ cf/2 & fb/4   \end{pmatrix} - t\cdot \begin{pmatrix} 1& 0 \\ 0 & 1  \end{pmatrix}\right)= t^2 - \tfrac{fb}{2}t+ \tfrac{(fb)^2}{16} - \tfrac{acf^2}{4}.\end{equation}

The discriminant of $\chi(t)$ is ${\cal D} = \tfrac{1}{4}\left(\tfrac{fb}{2}\right)^2 - \left( \tfrac{(fb)^2}{16} - \tfrac{acf^2}{4}\right) = \tfrac{acf^2}{4}$. We see that if $ac\ge 0$  (which is fulfilled by Proposition \ref{c2}) the discriminant is nonnegative  implying  the eigenvalues of $\tilde F_j^i$ are real.  Corollary \ref{jjj} is proved. 

\begin{Rem} \label{invar} 
For further use let us note that at the points such that  $ac=0$ the discriminant $\mathcal{D}$ of $\chi(t)$ given by  \eqref{2-1}   vanishes implying the   tensor $\tilde F_j^i$ has only one eigenvalue (of algebraic multiplicity two), namely $\tfrac{fb}{4}$.  At the points such that  $ac>0$ the discriminant $\mathcal{D}>0 $ implying the   tensor $\tilde F_j^i$ has two different real   eigenvalues. 
\end{Rem}

 \subsection{The function $L:= \sum_i \tilde F_i^i (:= \trace(\tilde F_j^i))$ is constant on each connected component  of the set of the points such that $B_1$  or $B_2$ is not defined. } \label{injo}

   \begin{Lemma}  \label{imp} 
   The function $L= \sum_i \tilde F_i^i (:= \trace(\tilde F_j^i))$ is constant on each connected component  of the set of the points such that $B_1$  is not defined. 
   \end{Lemma} 
   
   {\bf Proof.} Let at the point $p$ the form $B_1$ is not defined. We consider a small neighborhood $U(p)$ of $p$.  Lemma \ref{imp} is a direct corollary of the  following  
   
{  \bf Statement.  }  {\it  $L$ is constant on each connected component of the set $\{q\in U(p)\mid \textrm{ $B_1$ is not defined } \}  $.}  
   
   Now, the above  statement follows from the following two propositions: 
   
   \begin{Prop} \label{C1} Assume $B_1$ is not defined at every point of a  neighborhood of $p$. Then, $L$ is constant in this neighborhood. 
   \end{Prop} 
   
   \begin{Prop} \label{C2} Assume every neighborhood of $p$ has a point such that $B_1$ is  defined. 
   Then,  for a certain neighborhood $U(p)$ the function   $L$ is constant on the connected component of the set $\{ q \in U(p) \mid   \textrm{ $B_1$  is not defined} \}$ 
    containing the point $p$.  \end{Prop} 
    
   We will proceed as follows: we will first prove Proposition \ref{C1}. Then, we prove a technical Proposition  \ref{ppp}. Finally, we will use  Propositions \ref{C1}, \ref{ppp} in the proof of Proposition \ref{C2}.

{\bf Proof of Proposition \ref{C1}.} Our goal is to prove that $dL= 0$ at $p$. Without loss of generality, by Remark \ref{rem3},  we can assume $c\ne 0$ at the  point $p$.  
  We assume that $dL\ne 0$ at $p$, and find a contradiction.     

We denote by $W_0$ the connected component of the set 
$W:= \{ q \in T^2 \mid \textrm{ $B_1$ is not defined } \}$ containing $p$.  
We denote by   $\alpha :(-\infty, +\infty)\to T^2$ 
the integral curve of $V_2$ such that $\alpha (0)=p$. 
Since $W$ is invariant with respect to the flow of $V_2$, the curve $\alpha$ is a curve on $W_0$.
 Let us show that the curve $\alpha$ is periodic. 

In a small neighborhood of every point of the curve, we have $a\equiv 0$ implying 
$$L= \trace{(\tilde F^{i}_j)} = \trace\left(\begin{pmatrix}0  & b/2 \\ b/2 & c \end{pmatrix}\begin{pmatrix}0  & f/2 \\ f/2 & 0 \end{pmatrix}  \right)  =fb/2. $$ 
Then, the second equation of \eqref{sys}   implies that,   on $W$,  the function $L$ is invariant with respect to  the flow of 
$V_2$. Hence,  at every point of the curve  $\alpha$  we have $dL\ne 0$. 

Then, the connected component of the set $\{q \in T^2\mid L(q)= L(p)\}$ containing $p$ coincides  with the image of $\alpha$. Since $\{q \in T^2\mid L(q)= L(p)\}$  is compact, the image of $\alpha $ is compact implying the image of the curve is a closed circle.

The following cases are possible: 
\begin{itemize} 
\item[] {\bf Case (a):}  For every  $t\in \mathbb{R} $,  the form $B_2$ is  defined at the point $\alpha(t)$,
\item[] {\bf Case (b):}  There exists $t\in \mathbb{R} $ such that at the point $\alpha(t)$ the form $B_2$ is not defined. 
\end{itemize} 

Under assumptions of Case (a), 
let us construct a perfect coordinate system in a neighborhood $U(\alpha(t))$ of every point $\alpha(t)$.
 We assume that every  neighborhood  $U(\alpha(t))$ is sufficiently small  and  is homeomorphic to the disk.

As the first coordinate $x$ we take the function $-{2}L$  (where $L=  \sum_{i,j}\tilde F^{ij}g_{ij}$  as above).  Since $L$ is preserved by flow of $V_2$, its differential is not zero in a small neighborhood of every point  $\alpha(t)$.  Since 
$dL(V_2)=0$,  the coordinate $x$ can be taken as the first  admissible coordinate. 

In order to construct the second 
coordinate $y$, we consider the curve $\gamma:[0,t+1]\to W$  connecting the points $p= \gamma(0) $ 
and $q\in U(c(t))$  the such that $\gamma_{|[0,t]}= \alpha_{|[0,t]}$, $\gamma(t+1)=q$,  and such that $\gamma_{|[t, t+1]}$ lies in $U(\alpha(t))$.    We put   $y(q):= \int_{\gamma}B_2$. 

The function $y$  is well-defined, its differential is  $B_2$ and is not zero at $\alpha(t)$.  
The local  coordinates $x,y$ are  as in Proposition \ref{c4}(b). Then, in this coordinates, the metric $g$ is equal to  $  \left(X(x) {+\tfrac{1}{2}y}\right)dxdy$.  Since $V_2(X)= 0$ locally,  and since the functions $X(\alpha(t))$ coincides on  the intersection of the 
neighborhoods $U(\alpha(t_0))$ and $U(\alpha(t_0+ \varepsilon))$ (for small $\varepsilon$),  
 for every point of the curve $\alpha$ we have $X(\alpha(t))= X(\alpha(0))= X(p)$.

When $t$ ranges from $-\infty$ to $+\infty$, the coordinate 
$y$  also ranges from $-\infty $  to $+\infty$. Indeed,  $\int_{\alpha_{|{[0,t]}}}B_2(\alpha'(t))= \int_0^t {B_2}_{|\alpha(s)}({V_2}_{| \alpha(s)}) ds $, and $B_2(V_2) $ is positive and is therefore separated from zero on the compact set $\textrm{image}(\alpha)$. 

Then, there exists $t$ such that the value of $y$ corresponding to $\alpha(t)$ is $-2X(p)$ . At the  point $\alpha(t)$, the metric $g =  \left(X(x) {+\tfrac{1}{2}y}\right)dxdy$ is degenerate which contradicts the assumptions. Proposition \ref{C1} is proved under the additional assumptions of Case (a).

Let us now  prove Proposition \ref{C1}  under assumtions of Case (b): we 
assume that there exists  $t$ such $B_2=0$  at  $\alpha(t)$ this point .

Let $(t_{min}, t_{max})$, where $ t_{min} < 0 <  t_{max} \in \mathbb{R}$,   be the (open)  interval such that 
\begin{itemize} \item $B_2$ is defined at $\alpha(t)$ for  every $t\in (t_{min}, t_{max})$,  
 \item   $B_2$ is not defined at 
 $\alpha(t_{min})$, and at $\alpha(t_{max})$. \end{itemize}

As in the proof for Case (a), 
we  construct a perfect  local coordinate system $x,y$ in a neighborhood $U(\alpha(t))$ of every point $\alpha(t)$, where $t\in (t_{min}, t_{max})$. We put   $x(q):= -{2}L(q) $ and 
$y(q):= \int_{\gamma}B_2$, where   $\gamma:[0,t+1]\to W$, \ $\gamma_{|[0,t]}= \alpha_{|[0,t]},$ \ $\gamma(t+1)=q$,  and such that $\gamma_{|[t, t+1]}$ lies in $U(\alpha(t))$.  We assume that  the neighborhood $U(\alpha(t))$ is sufficiently small implying $B_2$ is defined at every point of $U(\alpha(t))$, and is homeomorphic to the disk. 
 
By Proposition \ref{c4}, in this coordinates, the metric is $(X(x) + \tfrac{1}{2}y)dxdy$. 
Let us show that the coordinate $y$  converges to $-2X(p)$ when $t$ converges to $t_{max}$. 

 In order to do this, we consider the scalar product on $T^*T^2$ induced by  $g$ (we will denote this scalar product by $g^*$). We consider the 
   function $h:= g^*(-{2}dL , B_2)$. This is indeed a function (i.e., $h$ does not depend on the choice of an  admissible  coordinate system) which is defined at the points such that $B_2$ is defined.  In admissible coordinates $(\tilde x, \tilde y)$ 
     in the neighborhood of the point $\alpha(t_{max})$, the function is given by 
   $ h= -{2}\tfrac{1}{\tilde f} \cdot \tfrac{\partial L }{\partial \tilde x}  \cdot  \tfrac{1}{\sqrt{|\tilde c|}}$. Since 
   $\tilde c(\alpha(t_{max}))=0$, we have  $h(\alpha(t))\stackrel{t \to t_{max}}{\longrightarrow} \pm \infty   $. 
   In the constructed above coordinates $(x,y)$, we have $h(\alpha(t))= \tfrac{1}{X(p)+ \tfrac{y(\alpha(t))}{2} }\cdot 1 \cdot   1$. Then,    ${X(p)}+ \tfrac{y(\alpha(t))}{2} \stackrel{t \to t_{max}}{\longrightarrow} 0$. Thus, $y(\alpha(t)) \stackrel{t \to t_{max}}{\longrightarrow}  -2X(p)$. 
   
  Similarly one can show that  the same is true for $t_{min}$, namely   $y(\alpha(t)) \stackrel{t \to t_{min}}{\longrightarrow}  -2X(p)$. . 
  
  Since $y(\alpha(t)) = \int_0^t {B_2}_{|\alpha(s)}({V_2}_{| \alpha(s)}) ds $,  and ${B_2}_{|\alpha(s)}({V_2}_{| \alpha(s)}) $ is positive for all $s\in (t_{min}, t_{max})$,  the values of $y(\alpha(t)) $  can not converge to the same number for $t\to t_{max}$ and for $t\to t_{min}$.  The  obtained  contradiction proves  Proposition \ref{C1}.

\begin{Prop} \label{ppp}  The set $\{q\in T^2 \mid \textrm{ $B_1$  or $B_2$ 
is defined in $q $  } \} $ is connected.  
\end{Prop}  

{\bf Proof.} It is sufficiently to prove that every point $p $ has a neighborhood $U(p)$ 
such that  the set  $S(p):= \{q\in U(p) \mid \textrm{ $B_1$  or $B_2$ 
is defined in $q $  } \} $ is connected.   We  take a sufficiently small  $U(p)$, and   consider admissible coordinates $x,y$ 
in  $U(p)$. We assume that the neighborhood is small enough so we can connect every two points  of this neighborhood by a geodesic. 

If  the set   $S(p)$ is not connected, at every point  $q\in U(p)$ we have 
$a(q)=0$, or $c(q)= 0$. Without loss of generality we can assume that at every point of $U(p)  $  we have $a=0$. Then, the point $p$  satisfies the assumptions of Proposition \ref{C1} above implying $L= fb/2= \const$  on $U(p)$.
 
 Let us now consider the points $U(p)\setminus S(p)$. At every such point, $a=c=0$ implying 
 $$ \tilde F^{ij} = \begin{pmatrix}0  & b/2 \\ b/2 & 0 \end{pmatrix} = \begin{pmatrix}0  & L/f \\ L/f & 0 \end{pmatrix} = \tfrac{{L}}{4} \cdot   g^{ij}.$$ Thus, at such points, $F = \tfrac{L}{2} H= \const \cdot H$. Without loss of generality we can assume that $\const =0$, otherwise we can replace $F$ by $(F- \const \cdot H)$. 
 
\begin{figure}[ht!] 
  %\hspace{-2cm}
  {\includegraphics[width=.3\textwidth]{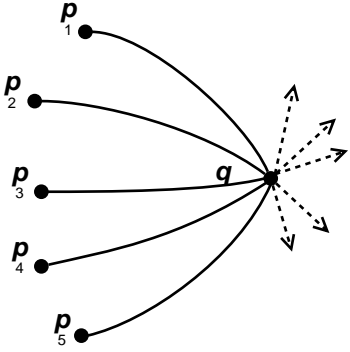}} 
 \caption{The  geodesic connecting the points $p_i$ with   the point $q$, and their  tangent vectors at the point $q$. For almost every $q$, the tangent vectors are mutually nonproportional}\label{fig}
\end{figure}

 We take 5 points  $p_1,...,p_5\in U(p)\setminus S(p)$ such that $F_{|T^*_{p_i}T^2}=0 $ at these points. 
  Since $F$ is an integral,  it  vanishes on every geodesic passing through any of the points $p_1,..., p_5$. Take a point $q\in U(p)$ in a small neighborhood of $S$, and   connect this point with  the points $p_1,...,p_5$   by geodesics, see Figure  \ref{fig}.  
 Let $ \xi_1\in T^*_{p_1}T^2, \dots  , \xi_5\in T^*_{p_5}T^2$ be the  vector-momenta  of these  geodesics at $q$.  At almost every $q$, the tangent vectors of the geodesics are mutually nonproportional implying    the   vector-momenta $\xi_i$ and $\xi_j$
  are not proportional for $i\ne j$. 
  
  Since $F$ is an integral and $F_{|T^*_{p_i}T^2}\equiv 0$, we have  $ F(\xi_i)= 0$. Thus,   the quadratic function $F_{|T_qT^2 }$ vanishes in 5 mutually nonproportional  points $\xi_i$. Hence,   $ F_{|T^*qT^2}\equiv  0$. Thus,  the restriction of $F$ to a small neighborhood of $p$ vanishes, which clearly contradicts the assumptions. The contradiction proves Proposition \ref{ppp}. 
 
 Combining Proposition \ref{ppp}, Remark \ref{4},   and Lemma \ref{impcomplex}, we obtain  
 \begin{Cor} \label{cc3} 
 Let $a>0$ at a point. Then, at every point of $T^2$ we have $a\ge 0$, $c\ge 0$. 
 \end{Cor}

{\bf Proof of Proposition \ref{C2}.}   We consider admissible coordinates $x,y$  in  a small neighborhood $U(p)$. We think that the point $p$ has the coordinates $(x(p), y(p))= (0,0)$. 
 In this coordinates, by Remark \ref{4}, the connected component of the  set  $\{q \in U(p) \mid \textrm{$B_1$ is not defined  at $q$}\} $  containing $p$  is one of the following sets (for a certain $\varepsilon>0$): 
 $$
 W_{+\varepsilon}:= \{q  \in U(p) \mid  0\le  x(q)\le \varepsilon \} \,,  \    W_{-\varepsilon}:= \{q  \in U(p) \mid  0\ge  x(q)\ge -\varepsilon \} \, , \  \textrm{or} \  W_{0}:= \{q  \in U(p) \mid   x(q)=0  \}. 
  $$ If the   connected component of the  set  $\{q \in U(p) \mid \textrm{$B_1$ is not defined  at $q$}\} $  containing $p$  is $W_{+\varepsilon}$ or $W_{-\varepsilon}$, we are done by  Proposition \ref{C1}. 
   We assume that the connected component of the  set  $\{q \in U(p) \mid \textrm{$B_1$ is not defined  at $q$}\} $  containing $p$  is $W_0$. Our goal is to prove that $\frac{\partial L}{\partial y}=0 $ for the points of this set.

 Let us first observe   that $da_{|q} = 0$ for every $q\in W_0$. 
 Indeed, by Corollary \ref{cc3}, the function $a$  accepts an extremum (minimum or maximum) at  $q$.

 Then, the second equation of \eqref{sys} tells us that $\frac{\partial L}{\partial y}=0$, i.e., $L$ is constant on the set $\{q\in U(p)\mid x(q)= x(p)\}$. Proposition \ref{C2} and  Lemma \ref{imp} are  proved.

\begin{Rem} \label{thesame}  Since there is no essential difference between $B_1$ and $B_2$, 
the function $L$ is constant on every connected component of the set  $\{ q\in {T}^2 \mid  {B_1}  \textrm{ or $B_2$  is not defined at  $q$}   \} $, as we claimed in the title of this section \end{Rem}

\subsection{ At  a neighborhood of every point the metrics are Liouville, or one eigenvalue of $\tilde F_j^i$  is constant on the manifold } \label{liou}

Recall that integrals linear in momenta and Killing vector fields  are closely related: 
 the function $I=\alpha(x,y) p_x + \beta(x,y) p_y$   is  an integral of the geodesic flow of $g$,  if and only if  
the vector field $v= (\alpha, \beta) $ is a Killing vector field. Moreover, the mapping 
$I=\alpha(x,y) p_x + \beta(x,y) p_y \mapsto v= (\alpha, \beta) $ is coordinate-independent.

 By Lemma \ref{impcomplex}, at every  point of $T^2$ we have $ac\ge  0$.  

\begin{Lemma} \label{twodiff} 
If there exists a point $q$ such that at this point at least one of the forms $B_1$, $B_2$ is not defined, then one of the eigenvalues of $\tilde F_j^i$ is constant on the manifold.   
\end{Lemma}  

\noindent{\bf    Proof. } We consider two sets:  
$$
W:= \{ p\in T^2 \mid \textrm{ $B_1$ and  $B_2$ are defined at $p$ }\} \textrm{ \ and \ } T^2 \setminus W. 
$$

Assume  $T^2 \setminus W \ne \varnothing$. At every point $s \in T^2$, we denote by $E_1(s)\le E_2(s)$  the roots of the characteristic polynomial 
$$
\chi(t):= \det(\tilde F_j^i - t \cdot \delta_j^i)  
$$
at the point $s$ counted with multiplicities.  (By Corollary \ref{jjj}, the roots of the polynomials $\chi(t)$ are real). 
 The functions $E_1$ and $E_2$ are at least continuous. 

At the points of $T^2 \setminus W$, by Remark \ref{invar}, we have $E_1=E_2= L/2$, where $L= \trace(\tilde F_j^i)$. 
Then, by Remark \ref{thesame},  both functions $E_1, E_2$ are constant on each connected component of  $T^2 \setminus W$. 

Since $W$ is open, and since $W\cup \left(T^2 \setminus W\right)= T^2$, in order to prove Lemma \ref{twodiff}, 
it is sufficient to show that  at least one of the functions $E_1, E_2$ is constant on every connected connected component of $W$.

We consider a point  $p$ such that at this point $ac>0$, and  denote by $W_0$  the connected 
 component of $W$ 
containing  $p$.

At every point $p_0$ of $W_0$, we consider the vector fields $\tfrac{\partial }{\partial x}$, $\tfrac{\partial }{\partial y}$, where $x,y$ are perfect coordinates is a neighborhood of $p_0$.     Though the perfect  coordinates are local coordinates, these vector fields are well defined at all points of $W_0$, see Remark \ref{pp}. Moreover, at every point $p_0$ the vectors  $\tfrac{\partial }{\partial x}$, $\tfrac{\partial }{\partial y}$  form a  dual basis  to the basis  $(B_1, B_2)$ in $T_{p_0}^*T^2$. 

Let us show that the vector fields $\tfrac{\partial }{\partial x}$, $\tfrac{\partial }{\partial y}$ are complete on $W_0$. Since the basis $(B_1, B_2)$ is  dual to the basis  $\left(\tfrac{\partial }{\partial x}, \tfrac{\partial }{\partial y}\right)$, it is sufficient to show that 
for every point $q$ of the boundary  $\partial W_0 := \overline  W_0 \setminus W_0$ 
the integral $\int_{p}^q B_1=\pm \infty$, or  $\int_{p}^q B_2=\pm  \infty$. 
We consider admissible coordinates $\tilde x, \tilde y$  in a neighborhood of $q$. 
Without loss of 
generality, $\tilde x(q)=  0$ and   $\tilde a(0)=  0$. 
As we explained in the proof of Lemma  \ref{imp},  the differential $d\tilde a_{|q}=0$ implying 
$\tilde a(\tilde x)=  \tilde x^2 \alpha(x)$, where  
$\alpha(x)$ is a smooth function in a neigborhood of $0$.   Then,  $\int_{p}^q B_1 = \const + \int_{\tilde x_0}^0 \tfrac{1}{\sqrt{|\tilde a(s)|}}  ds = \const \pm   \int_{\tilde x_0}^0\left(\tfrac{1}{{|s| \sqrt{|\alpha(s)|}}}  \right)  ds = \pm \infty$.

Thus, the vector fields $\tfrac{\partial }{\partial x}$, $\tfrac{\partial }{\partial y}$ are complete on $W_0$.

 We consider the local  coordinates  $u= \tfrac{1}{2}(x +y)$ and $v = \tfrac{1}{2}(x-y)$, and the corresponding vector fields 
 $\tfrac{\partial  }{\partial u}  =  \tfrac{1}{2}\left(  \tfrac{\partial  }{\partial x} +  \tfrac{\partial  }{\partial y}\right)$ and $\tfrac{\partial  }{\partial v}  =  \tfrac{1}{2}\left( \tfrac{\partial  }{\partial x} -  \tfrac{\partial  }{\partial y}\right)$. Since $\tfrac{ \partial }{\partial x}$,  $\tfrac{ \partial }{\partial y}$ are complete, the vector fields $\tfrac{\partial  }{\partial u}$, $\tfrac{\partial  }{\partial v}   $
are also complete.  

 The coordinates $u,v$ are as in Proposition \ref{c1}. Then, by Proposition \ref{c1},  in the coordinates $(u,v)$, the metric and the integral have the form $(U(u)- V(v))(du^2 - dv^2)$ and   $ \tfrac{U(u)p_v^2 -  V(v)p_{u}^2}{U(u)-V(v)}$. 
 Since $f= U(u)-V(v)>0$, we have $U(u)>V(v)$.  
 
 Let us note that  at every point of $W_0$, 
 the local functions $U$ and $V$ have a clear geometric sense, and, therefore, are globally given at all points of $W_0$, and can be continuously prolonged up to the boundary. 
 Indeed, in the coordinates $(u,v)$ the matrix of $\tilde F_j^i$ is  
 $$\begin{pmatrix}-V(v)  & 0 \\ 0 & -U(u) \end{pmatrix}. $$
Thus, $U= -E_1$ and $V= -E_2$.

   Consider the action of the group  
 $(\mathbb{R}^2, +)$ on $W_0$ generated by the vector fields $\tfrac{\partial  }{\partial u}   $  and $\tfrac{\partial  }{\partial v}$. The action  is well defined, since the vector fields commute and are complete.  The action is   transitive and locally-free. Then, $W_0$ is diffeomorphic to the torus, to the cylinder, or  to $\mathbb{R}^2$. 
  Since $T^2\setminus W_0 \ne \varnothing$, $W_0$ can not be the torus. 
 
Now suppose $W_0$ is a cylinder. Then, its boundary has at most two connected components. 
  Each integral curve of 
   at least one of the vector fields $\tfrac{\partial  }{\partial u}   $  and $\tfrac{\partial  }{\partial v}$ is  not closed. Without loss of generality, we 
   assume that  for every  $p\in W_0$ the integral curve of the vector field $\tfrac{\partial  }{\partial v}$ is not closed (i.e., it is the generator of the cylinder, or a standard winding on the cylinder. In the case the boundary of $W_0$ has two boundary components,   the integral curve of $\tfrac{\partial  }{\partial v}$ attracts to one component of the boundary for $t\to +\infty$, and to another component of the boundary for $t\to -\infty$). 
   
For every boundary component,    there exists a sequence of the points of any 
 integral curve of $\tfrac{\partial  }{\partial v}$ converging to a point of  the boundary component.  Indeed, the closure of $W_0$ is compact, so every sequence of points has a converging subsequence. We consider a  converging subsequence of the sequence 
  $\phi(0,p)=p, \phi(1,p), \phi(2,p),\phi(3,p),...$ where $\phi:\mathbb{R}\times W_0\to W_0$ 
 denotes the  flow of the vector field $\tfrac{\partial  }{\partial v}$. 
 Clearly, this  sequence  can  not 
  converge to a point of $W_0$. Then, it converges to a point of a boundary component.    
   Since the function $E_1=-U$ is constant along the integral curve,  the value of $E_1$ on the boundary coincides with the value of $E_1$ at the point $p$. Similarly, 
   the sequence 
 points $\phi(0,p)=p, \phi(-1,p), \phi(-2,p),\phi(-3,p),...$ has  a subsequence converging to another component of the boundary. Then, the value of $E_1$ on both components of the boundary coincides and is  equal to the value of $E_1$ at every point of $W_0$.   Then, the function $E_1$ is constant on $W_0$.

   Let us use the same idea to  show that $W_0$ can not be diffeomorphic to $\mathbb{ R}^2$. Indeed, in this  case   $\partial W_0$ has one connected component, and  the orbits of both vector fields $\tfrac{\partial  }{\partial u}$, $\tfrac{\partial  }{\partial v}$ are not closed implying $U(u)= V(v)$ at every point, which clearly contradicts the  assumptions.

   Finally, one of the eigenvalues of $\tilde F_j^i$ is constant on $W_0$.  Lemma \ref{twodiff} is proved. 
 
 \subsection{ If one eigenvalue of $\tilde F_j^i$ is constant, then there exists an integral linear in momenta} \label{eigenvalue}

 By Lemma \ref{twodiff}, we have the following two possibilities (not disjunkt): 
 
\begin{itemize} 
\item[{\bf (1) \ }] one of the eigenvalues of $\tilde F_j^i$ is constant,   
\item[{\bf (2) \ }]  at every point $ac>0$.   
\end{itemize} 

The goal of this section is to show that in the first case there exists an integral linear in momenta (at least on an appropriate  double cover of the torus; later (in \S \ref{4.2})
 we show that the integral  exists already on the torus, see Corollary \ref{cover}).  
 
 \begin{Lemma} \label{existlinear}  Let one of the eigenvalues of $\tilde F_j^i $ is constant. Then, for a certain (at most,    double) cover of the torus,  the lift of the integral is a linear combination of  the square of an integral linear in momenta and the lift of the Hamiltonian. Moreover, there exists no point $q$ such that $F_{|T^*_qT^2}\equiv \const\cdot  H_{|T^*_qT^2}$. 
 \end{Lemma} 
 {\bf Proof.} Without loss of generality we can assume that one of the eigenvalues of $\tilde F_j^i$ is identically $0$, otherwise we replace $F$ by $F- \const\cdot  H$ for the appropriate  $\const\in \mathbb{R}$.   Then,  
 $
 \tilde F_j^i
 $ has rank at most 1.

 Let $\tilde F_j^i \ne 0  $ at a point $q$. 
 We consider local  coordinate  $(u,v)$  in $U(q)$ such that  $\tfrac{\partial }{\partial u} $ lies in the kernel of $\tilde F_j^i$. In this coordinates, the (symmetric)  matrix of $\tilde F^{ij}$  satisfies the equation 
 $$
  \begin{pmatrix} \tilde F^{11} & \tilde F^{12}  \\  \tilde F^{21} &  \tilde F^{22}  \end{pmatrix} \begin{pmatrix} 1 \\ 0  \end{pmatrix}   =  0 
 $$ 
   implying $\tilde F^{11} = \tilde F^{12}= \tilde F^{21} =0$. 
   Then, in this coordinates $F= \tilde  F^{22} p_v^2 $ implying that the integral is locally the square of the function  
   $\sqrt{F^{22}} p_v$ , if $\tilde F^{22}>0$, or  $\sqrt{-F^{22}} p_v$,   if $\tilde F^{22}<0$. Then, the (linear in momenta) function  $\sqrt{F^{22}} p_v$ (if $\tilde F^{22}>0$) or $\sqrt{-F^{22}} p_v$    (if $\tilde F^{22}<0$) in  a local integral linear in momenta, and $\sqrt{F^{22}}\tfrac{\partial }{\partial v}$  (if $\tilde F^{22}>0$) or  $\sqrt{-F^{22}}\tfrac{\partial }{\partial v}$  (if $\tilde F^{22}<0$) is a Killing vector field. 
   
\begin{center}
\begin{figure}[ht!] 
  %\hspace{-2cm}
  {\includegraphics[width=.3\textwidth]{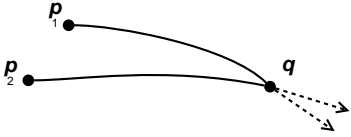}} 
  \caption{The  geodesic connecting the points $p_i$ with   the point $q$, and their  tangent vectors at the point $q$. At almost every $q$, the tangent vectors of the geodesics at the point $q$ are linearly independent}\label{fig3}
\end{figure}
\end{center}

   Let us show that the  points such that $\tilde F_j^i=0$ are  isolated. Indeed, otherwise there exist two such points, say $p_1$ and $p_2$,   in a sufficiently small neigborhood $U$. 
   For every   point $q$ of this neighborhood we  consider the geodesics connecting $p_i$ with q. 
   For almost every $q$, the geodesics intersect transversally at the point $q$, see Figure  \ref{fig3}.

  We denote by $\xi_1, \xi_2$ vector-momenta of these geodesics at the point $q$. 
   Since $F_{|T^*_{p_i}T^2}\equiv 0$, we have $F(\xi_1)=F(\xi_2)=0$ implyling
   $F_{|T^*_{p_i}T^2}\equiv 0$. Since this is fulfilled for almost every point $q$    of a small neighborhood,  the integral $F$ vanishes identically on two linearly independent vector-momenta,  which impossible for the integral $F= \tilde  F^{22} p_v^2 $ (for $\tilde  F^{22}\ne 0$). 
   
   Thus, the points $q$ such that ${ F}_{|T^*qT^2 }\equiv 0$  are isolated.
 Then, the set   $N:= \{ q\in T^2 \mid { F}_{|T^*qT^2 }\equiv 0\}$ is discrete.
       Hence, the set 
   $T^2\setminus N = \{ q\in T^2 \mid F_{|T^*qT^2 }\not \equiv 0\}$ is  connected implying that $\tilde F^{ij}$ is nonpositive definite  everywhere, or nonnegative definite   everywhere. Without loss of generality we can think that   $\tilde F^{ij}$ is nonnegative definite  everywhere, otherwise we replace $F$ by $-F$. 
   
   Let us show that in  {\it a small  neighborhood $U(p)$ of every point $p$ there exists precisely two integrals linear in momenta such 
   that \begin{itemize} \item[(a)] they  are smooth at every  points     $q \not\in  N$, and \item[(b)] the square of each of  these integrals is equal to $F$. \end{itemize}}

If $p\not\in N$, the statement is evident: in the constructed above local coordianates $u,v$ the integrals are    $\pm \sqrt{F}= \pm \sqrt{ \tilde  F^{22} p_v^2 }=  \pm p_v \sqrt{ \tilde  F^{22}}  .$ Since every neighborhood has a point from $T^2 \setminus N$, in a neighborhood of every point there exist at  
 most two such integrals. Thus, in order to prove the statement
 above  we need to prove that in a neighborhood of every point from $N$ there exists at least one  such integral (the second one will be minus the first).

\begin{center}
\begin{figure}[ht!] 
  %\hspace{-2cm}
  {\includegraphics[width=.7\textwidth]{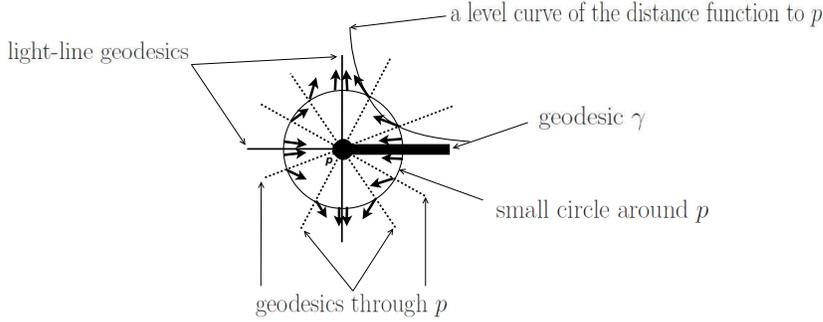}} 
  \caption{Qualitative behaviour of the vector field $v$ at the points  of  a small circle around $p$}\label{fig4}
\end{figure}
\end{center}

   Let $p\in N$. We take a small neighborhood $U(p)$ homeomorphic to the disk, and consider $U(p) \setminus \gamma$, where $\gamma$ is a geodesics starting at the point $p$, see Figure  \ref{fig4}. Since $U(p)\setminus \gamma$ is simply-connected and contains no point from $N$, on $U(p)\setminus \gamma$ there exists an integral $I= \alpha(x,y) p_x + \beta(x,y) p_y$   linear in momenta such that $I^2=F$.
   We consider the Killing vector field $v:= (\alpha, \beta)$ corresponding to  this integral. Since the value of this integral on each geodesic passing through $p$ is  zero, the Killing vector field $(\alpha, \beta)$ is orthogonal to geodesics containing $p$. Then, the qualitative behaviour of the vector field at the points  of  a small circle around $p$ is as on Figure  \ref{fig4}.  
   Indeed, they are tangent to the level curves of the geodesic distance function to the point $p$, which are hyperbolas (one of them is on Figure \ref{fig4}) and  light-line geodesics though $p$.   
   We see that the vector field $v$ is oriented in the same  direction on the different sides of  $\gamma$ , implying that one can prolong the vector field to $U(p)\setminus \{p\}$. Then,  there exists the integral $I$ linear in momenta  such that $I^2= F$ in $U(p)\setminus p$  as we claimed.  
   
   Since in 
  a small  neighborhood $U(p)$ of every point $p$ there exists precisely two integrals linear in momenta
   satisfying the conditions (a), (b) above,  an integral linear in momenta
   satisfying the conditions (a), (b) above exists on $T^2$,  or on the double cover of $T^2$.  The first statement of    Lemma \ref{existlinear} is proved.

  Let us prove the second statement of Lemma \ref{existlinear}:  let us show that the set $N$ is actually empty. Indeed, the index of the vector field 
   $v$ is  negative at  the points of $N$, see Figure  \ref{fig4}, and is zero at all other points. But the sum of the indexes of any vector field on the torus must be zero. 
   
   Thus,  there exists  an integral linear in momenta
   satisfying the condition (b) above on the torus, or on the double cover of the torus. 
   Lemma \ref{existlinear} is proved. 
    
    \begin{Cor} \label{zeros} 
    Let $v$ be a nontrivial  Killing vector field  of  a pseudo-Riemannian metric $g$ on the torus $T^2$. Then,      there is no point $p\in T^2$ such that $v=0$ at $p$.  
    \end{Cor}
     
     {\bf Proof.} In the Riemannian case (and, therefore, if $g$ has signature (--,--)), Corollary  \ref{zeros} is evident. Indeed, the Killing vector field preserves the complex structure corresponding to the metric, and is therefore holomorphic (with respect to the complex structure). By the Abel Lemma, it has  no zeros.   
    
    Let now  the signature of the metric be (+,--). We consider the integral linear in momenta  corresponding to the Killing vector field. It vanishes at the points where the Killing vector field vanishes. 
    The square of this  integral is an integral quadratic in momenta.  If the linear integral is 
    $    \alpha(x,y) p_x + \beta(x,y) p_y$, its square is $F= \alpha^2 p_x^2 + 2 \alpha \beta p_xp_y + \beta^2 p_y^2$, and the matrix $\tilde F^{ij}$  (such that $F= \sum_{i,j}\tilde F^{ij} p_ip_j$) is  
    $$
    \begin{pmatrix} 
    \alpha^2 & \alpha\beta \\ \alpha \beta & \beta^2 
    \end{pmatrix}. 
    $$
    We see that its rang is $\le 1$ implying that $0$ is a (constant) eigenvalue of $\tilde F_j^i$. 
    Then, by Lemma \ref{existlinear}, there exists no point such that $\alpha= \beta=0$ implying there exists no point such that  $v=  0$. Corollary  \ref{zeros}  is proved.

    \begin{Rem} Actually, our final goal is to prove that a nontrivial integral linear in momenta exists already on the torus (and not on the double cover of the torus). We will do it later, in Section \ref{proof4}. By Corollary \ref{cover} (whose proof does not use Theorem \ref{main1}, so no logical loop appears), the integral linear in momenta
   satisfying the condition (b) above exists already on the torus. \end{Rem}

\subsection{Proof of Theorem 2 under  the assumption that the vector fields $V_1$, $V_2$ exist  on the whole torus} \label{endof}

 Let the geodesic flow of  $g$ of signature (+,--) on the torus admits an integral quadratic in momenta; assume  the integral is not a  linear combination of the square of an integral linear in momenta and the Hamiltonian.  
 As everywhere in  Section \ref{4},  we assume that the vector fields $V_1$, $V_2$   satisfying conditions (A,B,C) from   \S \ref{admissible} exist on the whole torus.     By  Lemmas \ref{twodiff}, \ref{existlinear}, at every point of the manifold $ac >0$.

 We consider the vector fields $\tfrac{\partial  }{\partial u}$, $\tfrac{\partial  }{\partial v}$ from the proof of Lemma \ref{twodiff}.   These vector fields commute  and never vanish. Then, they generate a locally free  action of $(\mathbb{R}^2, +)$ on $T^2$. The stabilizer $G$ of this actions   is a subgroup of $(\mathbb{R}^2, +)$  with the following properties: it is 
  \begin{itemize} 
\item discrete, and 
\item the quotient space is compact.  
\end{itemize}  
 Then, it is a lattice, i.e., $G= \{ k\cdot  \xi  + m\cdot \eta \mid (k,m)\in \mathbb{R}\}$ for certain linearly independent vectors $\xi, \eta$. Then,  there exists a  natural diffeomorphism   $\phi: \mathbb{R}^2/G \to   T^2$.  We identify  $\mathbb{R}^2/G$ and 
 $   T^2$  by this diffeomorphism and consider the lift of the metric and 
  the integral to $\mathbb{R}^2$.  By Proposition \ref{c1}, in the  coordinate system $(u,v)$ 
   on $\mathbb{R}^2$, the metric and the integral are  $(U(u)- V(v))(du^2 - dv^2)$ and   $\pm  \tfrac{U(u)p_v^2 -  V(v)p_{u}^2}{U(u)-V(v)}$, i.e., are as in Model Example 1. 
   Since the metric and the integral are preserved by the lattice, the functions $U$ and $V$ are preserved by the lattice as well.  Thus, the metric on  $\mathbb{R}^2/G $ are as in Model Example 1. Theorem \ref{main1}  is proved (under  the  additional assumption that the vector fields $V_1$, $V_2$ exist  on the whole torus).

\section{ Proof of Theorem 4,  final step of the proof of Theorem 2, 
and  proof of Theorem 3}  \label{proof4}

\subsection{ Flat metrics of signature (+,--) on $T^2$, and their Killing vector fields}  \label{flat}
By the Gauss-Bonnet Theorem, a metric of constant curvature on the torus is flat (= has zero curvature). 
Recall that by the { \it standard flat torus } we  consider $(\mathbb{R}^2/G, dxdy)$, where $(x,y)$ are    the standard coordinates on $\mathbb{R}^2$, and $G$ is a lattice generated by two  linearly independent vectors.

It is well-known that every  torus $(T^2, g)$ such that the metric $g$  is flat and has signature (+,--) is isometric to a  standard one. Indeed, by \cite{carriere},  the flat torus is geodesically complete implying its universal cover is 
isometric to $(\mathbb{R}^2, dxdy)$. The fundamental group of the torus, 
$(\mathbb{Z}^2,+)$, acts on   $(\mathbb{R}^2, dxdy)$.  The action is isometric, free, and discrete. It is easy to 
 see that every   orientation-preserving  isometry of $(\mathbb{R}^2, dxdy) $ without fixed points is a translation. Then, $\mathbb{Z}^2$ acts  as a lattice generated by two  linearly independent vectors, and $(T^2, g)$ is isometric to a certain $(\mathbb{R}^2/G, dxdy)$.

 The space of Killing vector fields of    $(\mathbb{R}^2, dxdy)$ is a 3-dimensional  linear vector space generated by two translations $(1,0) = \tfrac{\partial }{\partial x} $ and $(0,1) = \tfrac{\partial }{\partial y} $, and the  pseudo-rotation $(y,x) =y\tfrac{\partial }{\partial x}+y\tfrac{\partial }{\partial y}$. Then, the space of Killing vector fields on the flat torus $(\mathbb{R}^2/G, dxdy)$ is two-dimensional and is generated by the Killing vector fields $(1,0) = \tfrac{\partial }{\partial x} $ and $(0,1) = \tfrac{\partial }{\partial y} $. Note that, depending on the values of the constants $(\const_1, \const_2)\ne (0,0)$,  every   integral curve of the Killing vector field $ \const_1\cdot \tfrac{\partial }{\partial x} + \const_2\cdot \tfrac{\partial }{\partial y} $ is either a closed curve, or 
 an everywhere dense winding on the torus.

\subsection{ Killing vector fields on the torus  of nonconstant curvature }   \label{4.2}

   \begin{Prop} \label{kil} 
   Let the metric $g$ of nonconstant curvature on the torus $T^2$ admits a nonzero Killing vector field $v$.  
   Then,   
    there  exists  a free action of the group $(\mathbb{R}/\mathbb{Z}, +)$ on the torus such that the infinitesimal generator of this action is proportional to the Killing vector $v$ with a constant coefficient of  proportionality. 
   \end{Prop}
 {\bf Proof.}  We denote by $R$ the scalar curvature  of $g$. 
By Corollary \ref{zeros},  the vector field $v$ has no zeros on $T^2$.  Then, the Killing vector field generates a locally-free  action of the group  $(\mathbb{R},+)$. 
 Let us prove that the Killing vector field (after an appropriate scaling) actually generates the action of the group 
  $SO_1=  \mathbb{R}/\mathbb{Z}$ without fixed points.

 Indeed, take a point $p$ such that $dR\ne 0$, and consider the orbit of the Killing vector field containing the point. Since the flow of a Killing vector field preserves  the curvature, at every point $q$ of the orbit  we have  $R(q)= R(p)$ and $dR\ne 0$.  Then, the orbit coincides with the connected component of the set $\{q \in T^2 \mid R(q)=R(p)\}$ containing  the point $p$ implying it is a circle. 
 
 We consider the action $\rho: \mathbb{R}\times T^2\to T^2$  of the  group $(\mathbb{R},+)$ generated by the flow of the vector field. 
 Since the orbit through $p$ is a circle, for certain $t_0>0$ we have $\rho(t_0, p)= p$ and 
 for no $t\in (0, t_0)$ $\rho(t, p)= p$.  Without loss of generality we can think that $t_0=1$, otherwise  we replace $v$ by  $t_0 \cdot v$.
 
 Since the action $\rho$  is isometric and orientation-preserving, it 
 commutes with the exponential mapping $\exp:TT^2\to T^2$. Then,  for every point $q\in T^2$ we have $\rho(1,q)= q$ and  $\rho(t,q)\ne  q$ for $t\in (0, 1)$. Thus, the action of the group  $(\mathbb{R}/\mathbb{Z}, +)$  is well-defined, and has no fixed points.    Proposition \ref{kil} is proved. 
 
 \begin{Cor}  \label{cover1} Let $v$ be a nonzero  Killing vector field on the torus $(T^2,g)$, where $g$ has signature $(+,-)$.  Then, there exists 
  no  involution $\sigma:T^2\to T^2$ without fixed point that preserves the orientation  and  the metric, and sends the vector field  $v$ to $-v$.  
 \end{Cor} 
 {\bf Proof.} If the metric $g$ has constant curvature,  as we have recalled  in \S \ref{flat},   the torus is isometric to $(\mathbb{R}^2/G, dxdy)$ for a lattice $G$ generated by two  linearly independent vectors $\xi$ and $\eta$, and the Killing vector field is 
 $\const_1 \cdot \xi +  \const_2 \cdot \eta$ for $(\const_1, \const_2)\ne (0,0)$.  The involution $\sigma$ without 
 fixed points  that preserves the orientation  and  the metric induces an isometry of $(\mathbb{R}^2, dxdy)$ without fixed points that preserves the orientation  and  the metric.  Such isometry is a translation and can not send the  Killing vector field     $\const_1 \cdot \xi +  \const_2 \cdot \eta$ to $-\left(\const_1 \cdot \xi +  \const_2 \cdot \eta\right)$. Corollary \ref{cover1} is proved under the
 assumption that $g$ has constant curvature.  
 
 Assume now that the curvature of $g$ is not constant. Then, by Proposition \ref{kil}, the Killing vector field (after the appropriate scaling) generates a free  action of $(\mathbb{R}/\mathbb{Z}, +)$ on $T^2$. 
 We consider the quotient space $T^2/_{(\mathbb{R}/\mathbb{Z})}.$ Since the action of $\mathbb{R}/\mathbb{Z}$ on $T^2$ is free, 
 the quotient space is a 1-dimensional closed manifold, i.e., is diffeomorphic to $S^1$. The orientation of the torus induces the orientation on $S^1$.  
 
   The involution $\sigma$ of the torus preserves the action,  the   orientation, and sends $v$ to $-v$. 
   Then, it inverses the orientation of $S^1 = T^2/_{(\mathbb{R}/\mathbb{Z})}.$  Then, it has a fix  point. We consider the orbit of  $\mathbb{R}/\mathbb{Z}$ corresponding to this point. The involution  $\sigma$  preserves this orbits and 
    changes the direction of the vector field $v$ on this orbit. Then, it has a fixed point  which contradicts the assumptions. The contradiction proves Corollary \ref{cover1}.  
 
  \begin{Cor}  \label{cover} Let $F$ be a nontrivial  integral quadratic in momenta for the geodesic flow of the metric $g$  on the torus $T^2$ and $\pi:\widetilde T^2 \to T^2   $ be a double cover of $T^2$. 
Assume  the lift of the integral to $\widetilde T^2$ is  a   linear combination of the square of  a function  linear in momenta and the lift of the  Hamiltonian. Then, the integral $F$ 
is  a   linear combination of the square of  an integral  linear in momenta and the  Hamiltonian
 \end{Cor} 
 
 {\bf Proof.} We consider the involution $\sigma:\widetilde T^2 \to \widetilde T^2$ corresponding to the cover:
 $\sigma(\tilde p)= \tilde q$    if $\pi(\tilde p)=\pi(\tilde q)$ and $\tilde p\ne \tilde q$. The involution preserves the lift of the Hamiltonian and of the integral. 

We consider the function  $I:T^*\widetilde T^2\to \mathbb{R}$ linear in momenta such that 
$ F=  \const_1\cdot H + \const_2\cdot I^2$, where $H$ and $F$ denote the lift of the Hamiltonian and the integral. Since the integral $F$ is nontrivial, $\const_2\ne 0$ implying $I$ is a nontrivial  integral (linear in momenta).    We consider the Killing vector field $v $  corresponding to the integral. Since the involution $ \sigma$ preserves  $H$ and $F$, it preserves $I^2= \tfrac{1}{\const_2}(F- \const_1\cdot H)$. Since by Proposition \ref{kil}  the vector field $v$ vanishes at no point, 
either $d\sigma(v)= v $ for all points, or $d\sigma(v)= -v$ for all points.  The second possibility is forbidden by Corollary \ref{cover1}. Then, $d\sigma(v)= v $ implying the integral $I$ on $\widetilde T^2$ induces an  integral $I$ (linear in momenta) on $T^2= \widetilde T^2/\sigma$ such that, on $T^2$, $ F=  \const_1\cdot H + \const_2\cdot I^2$.    Corollary \ref{cover} is proved.

 \subsection{ Proof of Theorem \ref{main4} } Let $F$ be an  integral linear in momenta  of the geodesic flow of a metric $g$  on the torus $T^2$. We denote by $v$ the corresponding Killing vector field.   
 We consider the action  $\rho$  of  
 $(\mathbb{R}/\mathbb{Z}, +)$ on $T^2$ from  Proposition \ref{kil},   the   quotient space $T^2/_{(\mathbb{R}/\mathbb{Z})}$ diffeomorphic to the circle, and the  tautological   projection $\pi:T^2 \to T^2/_{(\mathbb{R}/\mathbb{Z})}=S^1 $. 
 Let us construct a coordinate system $(x  \in \mathbb{R }  \textrm{\ mod} \  1, \ y    \in \mathbb{R }  \textrm{\ mod } \  1 )$  on $T^2$. We parametrize   $S^1$  by   $(Y\in \mathbb{R}    \ \textrm{ \ mod \ } \  1 )$, and 
  put $y(q):= Y(\pi(q))\in  \mathbb{R}/\mathbb{Z})$. In order to construct the coordinate $x$, we consider a smooth section $c: S^1 \to T^2$   of the bundle. By definition of the section, for every $q\in T^2$ there exists  a unique $t\in (\mathbb{R }  \textrm{\ mod\ } \  1)$ such that $\rho(t,q) \in \textrm{image}(c)$.  We put $x(q)= -t$.   
 
 By construction, in this coordinates, the vector field $v$  is $\tfrac{\partial }{\partial x}$, and the corresponding integral linear in momenta is $p_x$.   Let in this coordiantes the metric $g$ is given by 
 $g= K(x,y) dx^2+ 2 L(x,y) dxdy + M(x,y) dy^2$. Since the metric  has signature (+,--), we have   $KM-L^2 = \det\begin{pmatrix} K  &L\\ L & M\end{pmatrix} <0$.  Thus, in order to prove Theorem \ref{main4}, it is sufficient to show that the functions $K, L, M$ are functions of the variable $y$ only, i.e., $ \tfrac{\partial K}{\partial x} =  \tfrac{\partial L}{\partial x}= \tfrac{\partial M}{\partial x}=0.$
 
 We denote by $k(x,y), l(x,y), m(x,y) $ the components of the inverse matrix to $g$: 
 $$
 \begin{pmatrix} k  &l\\ l &m\end{pmatrix} = \begin{pmatrix} K  &L\\ L & M\end{pmatrix} ^{-1}. 
 $$
 Evidently, $2H=  k(x,y) p_x^2 + 2 l(x,y) p_xp_y + m(x,y) p_y^2$, and 
  the condition $\{F,2H \}=0$ reads \\
 \begin{eqnarray*}
0&=& \left\{ p_x, k(x,y) p_x^2 + 2 l(x,y) p_xp_y + m(x,y) p_y^2\right\}  \\
 &=&  \tfrac{\partial k}{\partial x} p_x^2 +  2 \tfrac{\partial l}{\partial x} p_xp_y + \tfrac{\partial m}{\partial x} p_y^2 \, ,
 \end{eqnarray*}
 i.e., is equivalent to the condition $ \tfrac{\partial k}{\partial x} =  \tfrac{\partial l}{\partial x}= \tfrac{\partial m}{\partial x}=0.$ Then, the coefficients $k,l,m$ depend on the variable $y$ only, implying that the coefficients $K,L, M$  also depend on the variable $y$ only.   Theorem \ref{main4} is proved.

  \subsection{ Proof of Theorem \ref{main1} under the assumption that the vector fields $V_1, V_2$ do not exist on the torus }  \label{laststep} 
 We assume that the geodesic flow of the metric $g$ on $T^2$ admits a nontrivial integral $F$ quadratic in momenta that is not a linear combination of the Hamiltonian and an integral linear in momenta.    Assume  the vector fields $V_1,V_2$ satisfying assumptions (A,B,C) from \S \ref{admissible}   do not exist. We consider   the double cover $\pi: \widetilde T^2 \to T^2$ such that $V_1,V_2$ satisfying (A,B,C) exist on $\widetilde T^2$. Then, by the proved
  part of Theorem \ref{main1}, the lift of the  metric to $\widetilde T^2$ is as in  Model Example 1
  (we idientify  
   $\widetilde T^2$ with $\mathbb{R}^2/G$  and the lift $\tilde g$  of the 
    metric with the metric from Model Example 1).    On the torus $\widetilde T^2$, the  only  possibility for the vector fields $V_1,V_2$ are (we consider the standard orientation on $\mathbb{R}^2$): 
    $$ V_2 = \lambda \left(\tfrac{\partial }{\partial x} + \tfrac{\partial }{\partial y}\right), \ V_1 = \mu\left(\tfrac{\partial }{\partial x} - \tfrac{\partial }{\partial y}\right), $$
    where $\lambda$ and $\mu$ are smooth  functions on $\widetilde T^2$ such that for every $\tilde p\in \widetilde T^2$ we have  $\lambda(\tilde p) \mu(\tilde p)>0$, and $x,y$  are the standard coordinates on $\mathbb{R}^2$. 
    
    We consider the involution $\sigma$ 
    corresponding to the cover $\pi$,  that it  $\sigma(\tilde p)= \tilde q$ if and only if $\pi(\tilde p)= \pi(\tilde q)$ and  $\tilde p\ne \tilde q$.  Since by assumptions the vector fields $V_1, V_2$ 
     do not exist on $T^2$, and the involution  preserves the  orientation, the metric $\tilde g$,   and the lift of the integral, we have 
     $$
     d\sigma\left(\tfrac{\partial }{\partial x} + \tfrac{\partial }{\partial y}\right) = -\left(\tfrac{\partial }{\partial x} + \tfrac{\partial }{\partial y}\right) \ \ \textrm{and} \ \   d\sigma\left(\tfrac{\partial }{\partial x} - \tfrac{\partial }{\partial y}\right) = -\left(\tfrac{\partial }{\partial x} - \tfrac{\partial }{\partial y}\right) $$
   implying  \begin{equation}\label{inversion}
     d\sigma\left(\tfrac{\partial }{\partial x} \right) = -\tfrac{\partial }{\partial x} \ \ \textrm{and} \ \   d\sigma\left( \tfrac{\partial }{\partial y}\right) = -\tfrac{\partial }{\partial y}.  \end{equation} 
But on the torus $\mathbb{R}^2/G$ there is no involution with no fixed point with the property  \eqref{inversion}. The contradiction shows  that the situation assumed in this section, namely that the  vector fields $V_1, V_2$ do not exist on  $T^2$, is impossible.  Theorem \ref{main1}
is proved.

\subsection{ Proof of Theorem  \ref{main3}} 

We assume that $g$ is a metric of signature (+,--)  on the Klein bottle $K^2$ whose geodesic flow admits an integral quadratic in momenta. We also assume that  the lift  of the integral   
to  the oriented cover  is not a linear combination of the lift of the Hamiltonian and the square of a function linear in momenta.   Our goal is to prove that $(K^2, g)$   is as in Model Example 2.

We consider the oriented cover $\pi:T^2 \to K^2$, and   the lift of the metric and the integral to $T^2$. They satisfy the assumptions in Theorem \ref{main1}. Hence  we can think that $T^2$, the lift of the metric, and the lift of the integral are as Model Example 1: 
$$T^2= \mathbb{R}^2/G\, , \  g=(X(x)- Y(y))(dx^2 - dy^2)  \, , \  \textrm{and} \ 
F= \tfrac{X(x)p_y^2 -  Y(y)p_{x}^2}{X(x)-Y(y)}, $$ 
where $G= \{k\cdot \xi + m \cdot \eta\mid k,m\in \mathbb{Z} \}$. 
 
Next, consider the universal cover $\tilde \pi :=  \pi\circ P: \mathbb{R}^2 \to K^2$, where   $P$ is the canonical projection from $\mathbb{R}^2$   to $\mathbb{R}^2/G$. 
We conisder  the action of the fundamental group of the Klein bottle on $\mathbb{R}^2$ corresponding to $\tilde \pi$. Recall that the fundamental group of $K^2$ is  generated by two elements, say $A$ and $B$, satisfying the relation $ABA^{-1}B={\bf 1}:$
\begin{equation} \label{pi}
\pi_1(K^2)= \<A,B | ABA^{-1}B={\bf 1}\> .
 \end{equation} 
 
  This action has the following properties: 
 
 \begin{itemize} 
 \item[(a)] It preserves the metric  and the integral, 
 \item[(b)] It is free and discrete. 
 \end{itemize} 
Let us show that the condition (a) implies the condition 
  \begin{itemize} 
 \item[($\textrm{a}'$)] For every element $\alpha \in  \pi_1(K^2)$ we have \begin{equation} d\alpha(\tfrac{\partial}{\partial y}) = \pm \tfrac{\partial}{\partial y} \, , \ \ d\alpha (\tfrac{\partial}{\partial x}) = \tfrac{\partial}{\partial x}. \label{nn} \end{equation}  
\end{itemize} 
 
 Indeed, 
since at every point $(x,y)\in \mathbb{R}^2$  the factor $X(x)-Y(y)\ne 0$, and since every nonempty  level $\{X= \const_1\}$  intersects with every   nonempty  level $\{Y= \const_2\}$, without loss of generality we can think that $X(x)>Y(y)$ for all $(x,y) \in \mathbb{R}^2$. 

Now, in the coordinates $x,y$, the matrix of $F^i_j$ is $\begin{pmatrix} - Y(y) &  \\ & - X(x)  \end{pmatrix}, $ so $\tilde F_{j}^i$ has eigenvalues 
$-X(x)$, $-Y(y)$.     

Since the action  preserves the metric and the integral, it preserves the eigenvalues $X,Y$ and the eigenspaces $\textrm{span}\left(\tfrac{\partial}{\partial y}\right)$ and $\textrm{span}\left(\tfrac{\partial}{\partial x}\right)$   of this eigenspaces. 
Since $g\left(\tfrac{\partial}{\partial x}, \tfrac{\partial}{\partial x}\right) = X(x)-Y(y),  $ and $\alpha$ preserves $X$ and $Y$, we have that 
 $g\left(d\alpha\left(\tfrac{\partial}{\partial x}\right), d\alpha\left(\tfrac{\partial}{\partial x}\right)\right) =  X(x)-Y(y)$ implying $ d\alpha\left(\tfrac{\partial}{\partial x}\right)= \pm \tfrac{\partial}{\partial x}$. The proof that  $ d\alpha\left(\tfrac{\partial}{\partial y}\right)= \pm \tfrac{\partial}{\partial y}$ is similar.

Thus, the action preserves the standard flat metric $dx^2 + dy^2$ on $\mathbb{R}^2$.  Then, the fundamental group of $K^2$ as a  crystallographic group. From the classification of   crystallographic groups \cite[\S1.7]{berger}, it follows that every  action of the group \eqref{pi} on $\mathbb{R}^2$ satisfying ($\textrm{a}'$, b) is 
generated by $A$, $A(x,y)= (x+ c, -y)$  and $B$, $B(x,y)=(x, y + d)$ for certain $c\ne 0 \ne d$, i.e., is an the Model Example 2. Theorem   \ref{main3} is proved.

{\bf Acknowledgement.}   The author  thanks
Deutsche Forschungsgemeinschaft
(Priority Program 1154 --- Global Differential Geometry and Research Training Group   1523 --- Quantum and Gravitational Fields)   and FSU Jena  for partial financial support, and D. Alekseevsky, O. Bauer, A. Bolsinov, G. Manno,  P. Mounoud,   G. Pucacco, and A. Zeghib  for useful discussions.

 \end{document}